\documentclass[a4paper, 12pt]{article}
\input{./preamble.tex}

\title{\textsc{\ta-tilting theory\,: a self-contained introduction}}
\author{Arne Johannsmann}

\begin{document}
\bibliographystyle{alpha}
	
	\maketitle
    \vspace{10em}
    \begin{center}

    \[
        \begin{tikzcd}
            & \indd{110}{011}{010} \ar[r] \ar[dr] 
            & \ind{110}{010} \ar[r] \ar[dr]
            & \ind{110}{100}  \ar[r]
            & 100 \ar[dr]
            & \\		
            \indd{110}{011}{001} \ar[ur] \ar[r] \ar[dr] 
            & \indd{110}{100}{001} \ar[urr, crossing over] \ar[dr] 
            & \ind{011}{010} \ar[r]
            & 010  \ar[rr]
            &
            & 0 \\
            & \ind{011}{001} \ar[ur, crossing over] \ar[rr, bend right] 
            & \ind{100}{001} \ar[r] \ar[uurr, bend right, crossing over]
            & 001 \ar[urr]
            &
            & \\
        \end{tikzcd}
	\]
    \vspace{8em}

        \textsf{Universität Bielefeld}

        \textsf{Fakultät für Mathematik}
    \end{center}
	
    \newpage
	\begin{abstract}
        The following is a self-contained introduction to the basics of \ta-tilting theory. We assume that the reader is familiar with Auslander-Reiten theory, but circumvent the need for the Brenner-Butler tilting theorem completely.
	\end{abstract}
\setcounter{tocdepth}{2}
\tableofcontents

\newpage

\subsection{Introduction}

This text was originally written as a master thesis, following up on a sequence of lectures with the title and topic of \emph{Noncommutative Algebra}, held by Prof. Dr. William\footnote{Bill.} Crawley-Boevey in Bielefeld University, 2019-2021. My thesis was originally handed in on the 11th July, 2022.

A few changes were made to the original text, primarily to clean it up, and to split the final section into its own paper. It largely contained a few computational results, which I felt the need to include, but which didn't contribute to the goal of introducing \ta-tilting theory. These computations will be published separately.

This is the second version of this text. Several arguments in the first version were incomplete, or wrong. Most importantly, the proof that \ta-rigid modules are determined by their $g$-vectors (\cref{4gVectorsDetermineTauRigid}) had a gap, and the proof that tilting modules have exactly $\# A$ summands (\cref{4partialTiltingBound}) relied on the false claim that the indecomposable projective modules form a basis of $K_0(A)$. Both proofs have been corrected, and neither relies on the Brenner-Butler theorem. Further, the equality $\# \SP(\ST) = \supportrank \ST$ now explicitly requires $\ST$ to be functorially finite, and the first version wrongly claimed that the torsion functor $t(-)$ is always left exact (\cref{2TorsionFunctorsNotExact}). Several proofs now include steps which were previously omitted, most notably the proof of \cref{4PTIsTiltingOverAAnn}, where the injectivity of the envelope is now proved. Finally, a number of smaller errors have been fixed.

The primary aim of this paper is to provide an accessible, and self-contained introduction to \ta-tilting theory. We do this, assuming only that the reader is familiar with basic Auslander-Reiten theory, as it might be taught in any series of lectures on the representation theory of finite dimensional $K$-algebras, quivers, and path algebras.

We have done our best to ensure that the following text can be used for self-studies, or to teach a course on \ta-tilting theory. Despite that, errors might have slipped in. Feel free to contact me if you find any, if you have questions, or if you got any use out of this text whatsoever.

The paper \emph{\ta-Tilting Theory} \cite{zbMATH06293659}, released in 2013 by Adachi, Iyama, and Reiten, can—true to its title—be considered the origin of modern \ta-tilting theory. The authors develop the foundation of \ta-tilting theory, and—in their own words—'complete' classical tilting theory from the viewpoint of mutation.

The subject of interest are the so-called \emph{support \ta-tilting modules}. We give a preliminary definition:
\begin{definition}
    Given a finite-dimensional $K$-algebra $A$ over an algebraically closed field $K$, consider a basic finite-dimensional left $A$-module $T$. We say that $T$ is \ta-rigid if 
    \[
        \Hom_A(T, \tau T) = 0.
    \]
    We say that $T$ is a \emph{support \ta-tilting module}, if it is \ta-rigid, and its number of summands equals its support-rank, i.e. the number of isomorphism classes of simple modules showing up in its composition series:
    \[
        \# T \overset{!}{=} \supportrank T.
    \]
\end{definition}

Why anyone would be interested in modules with these properties is not obvious. To provide some motivation: It turns out that the set of isomorphism classes of support \ta-tilting modules $\sttilt A$ of $A$ is in natural bijection to several other classes of interesting objects. This includes certain subcategories with desirable properties, such as the so-called torsion classes, which we will introduce in the first section.

The incredibly rich combinatorial and homological structure of support \ta-tilting modules then ensures that $\sttilt A$ possesses the structure of a regular(!) graph, whose edges can be interpreted as 'mutation', or as inclusion between the aforementioned torsion classes associated with support \ta-tilting modules.

We will introduce \ta-rigid modules, and explore their connection to functorially finite torsion classes in the second section. In the third section we borrow some tools from classical tilting theory, which we then use in the fourth section to develop the theory of support \ta-tilting modules.

The fifth and sixth section are devoted to mutation, and to a short look into current topics of research.

Our treatment of the introductory material follows Bill's notes \cite{NA2} quite closely. We recommend them to the interested reader. Starting from Section 4 we primarily follow \cite{zbMATH06293659}, and rework several proofs for our convenience, and to meet our goal of making this exposition as self-contained as possible.

I am dedicating this thesis to my parents.

\subsection{Notation and Preliminaries}
Throughout this work we will consider finite-dimensional $K$-algebras (often denoted $A$), where $K$ is an algebraically closed field. When we are talking about $A$-modules, then we are referring to objects of $A\Mod$, the category of finite-dimensional left $A$-modules.

If $A$ fulfills that $A / \rad A$ is isomorphic to $K^n$ for some value of $n$, then it is possible to find some quiver $Q$, and an admissible ideal $I \subset KQ^+$, such that $A \cong KQ/I$. If it does not, then its module category $A\Mod$ is at least equivalent to the module category of such an algebra. For this reason we generally think of algebras in terms of their underlying quivers, and denote their modules by their dimension vectors, assuming this is unambiguous. For clarity we underline the projectives, and overline the injectives, for example, we might write $\ol{\ul{111}}$ for the projective-injective indecomposable of the path algebra of the linearly oriented $\mathbb{A}_3$ quiver.

\begin{definition}
    A \emph{module class} $\SC$ is a full subcategory of a module category, such that 
    \begin{enumerate}
        \ii $\SC$ is closed under isomorphisms
        \ii $X, Y \in \SC$ holds if and only if $X \oplus Y \in \SC$
\end{enumerate}
    Given a module $X$, we write $\add X$ for the smallest module class containing $X$. In particular, this is closed under taking direct sums. Every subcategory of $A\Mod$ which we encounter will be a module class. This primarily includes $\gen X$, the class of modules which can be written as a quotient of a module in $\add X$, and 
    \[
        \bong \SC := \{ Y \in A\Mod \mid \Hom_A(Y, \SC) = 0 \}.
    \]
    Each of these will be introduced as it becomes necessary.
\end{definition}

Given some $A$-module $X$, we will write $\# X$ for the number of isomorphism-classes of indecomposable summands of $X$, that is, the number of indecomposable summands counted without multiplicities. $\# A$ is the number of vertices of $A$.

We expect the reader to be familiar with Auslander-Reiten theory:
\begin{theorem*}[Auslander-Reiten Formula]
    Given $A$-modules $X, Y$, we obtain the following isomorphisms:
    \[
        \ul{\Hom}_A(\tau^-Y, X) \cong D\Ext_A^1(X, Y) \cong \ol{\Hom}_A(Y, \tau X).
    \]
\end{theorem*}

In particular, we quickly mention the following result, which we use thrice.
\begin{proposition}\label{1ProjDimEquiv}
    Let $M$ be an $A$-module. Then the following statements are equivalent:
    \begin{enumerate}
        \ii $\projdim M \leq 1$
        \ii $\Hom_A(DA, \tau M) = 0$
        \ii For all injective $A$-modules $I$, we have that $\Hom_A(I, \tau M) = 0$
    \end{enumerate}
    \begin{proof}
        Take a minimal projective presentation $P_1 \xrightarrow{f} P_0 \twoheadrightarrow M$. Now, $\projdim M \leq 1$ if and only if $\ker f = 0$. We have an exact sequence 
        \[
            \begin{tikzcd}
                0 \rar & \tau M \rar & \nu P_1 \rar["\nu f"] & \nu P_0
            \end{tikzcd}
        \]
        Applying $\nu^- = \Hom_A(DA, -)$ yields a commutative diagram with exact rows.
        \[
            \begin{tikzcd}
                0 \rar & \Hom_A(DA, \tau M) \rar & \nu^- \nu P_1 \dar["\cong"] \rar["\nu^- \nu f"] & \nu^- \nu P_0 \dar["\cong"] \\
                0 \rar & \ker f \rar & P_1 \rar["f"] & P_0 
            \end{tikzcd}
        \]
        We obtain that $\ker f = \Hom_A(DA, \tau M)$, and the statement follows.
    \end{proof}
\end{proposition}

\section{Torsion Theories}

We introduce torsion theories, and discuss their fundamental properties. As we will later see, certain classes of particularly nice torsion theories share a close and deep connection to support \ta-tilting modules, making this topic one of the building blocks of the theory.
\subsection{Definition and First Properties}

\begin{definition}\label{2TorsionTheory}
	A \emph{torsion theory} in a module category $A\Mod$ is a pair of module classes $(\ST, \SF)$—the \emph{torsion} and \emph{torsion-free} classes—such that:
	\begin{enumerate}
		\ii $\Hom_A(\ST, \SF) = 0$
		\ii Any object $X \in A\Mod$ fits into an exact sequence of the form
		 \[
			\begin{tikzcd}
				0 \rar & tX \rar & X \rar & X/tX \rar & 0
			\end{tikzcd}
		\] 
		where $tX \in \ST,~X/tX \in \SF$.
	\end{enumerate}
    This definition—like many others—extends to general abelian categories, but for our purposes we only need to think in terms of module categories.
\end{definition}

\begin{remark}
    It is easy to see that:
    \begin{enumerate}
        \ii Take modules $T \in \ST,~X \in A\Mod,~F \in \SF$. Any map $T \rightarrow X$ factors through $tX$. Dually, any map $X \rightarrow F$ factors through $X/tX$.
        \ii Given some torsion theory $(\ST, \SF)$ any simple module $S$ is either contained in $\ST$, or in $\SF$.
    \end{enumerate}
\end{remark}

\begin{notation*}
	For a set $\SC$ of modules in $A\Mod$ we define:
	\begin{enumerate}
		\item $\SC^{\perp i,j,\dots} := \{ X \mid \Ext^{n}_A(M,X) = 0 \text{ for all } M \in \SC \text{ and } n = i, j , \dots\}$
		\item ${}^{\perp i,j,\dots}\SC := \{ X \mid \Ext^{n}_A(X,M) = 0 \text{ for all } M \in \SC \text{ and } n = i, j , \dots\}$
	\end{enumerate}
	Of course, $\Ext^0_A = \Hom_A$. In practice our indices will only take on the values $0$ or $1$, and thus we will never need more than two at a time.

    If $\SC$ contains a single module $M$, or can be written as $\add M$ for some module $M$, then we will write $M^{\perp i,j,\dots}$ or ${}^{\perp i, j,\dots}M$ to refer to the above. Since $\Ext$ preserves finite direct sums in both arguments, it is clear that $M^{\perp i,j,\dots} = (\add M)^{\perp i,j,\dots}$, as well as that ${}^{\perp i, j,\dots}M = {}^{\perp i, j,\dots}\add M$. The same property shows that $\SC^{\perp i,j,\dots}$, and ${}^{\perp i,j,\dots}\SC$ are module classes.
\end{notation*}

\begin{proposition}\label{2Adjoints}
	Let $(\ST, \SF)$ be a torsion theory. 
	\begin{enumerate}
		\ii $\ST = {}^{\perp 0} \SF$, as well as $\SF = \ST^{\perp 0}$. More precisely, there are inclusion-reversing bijections \[ \left\{
		\begin{array}{c}
			\text{Torsion} \\
			\text{classes}
		\end{array} \right\}
		\begin{array}{c}
			\xrightarrow{\enspace \, \ST^{\perp 0} \,\enspace} \\
			\xleftarrow{\enspace {}^{\perp 0}\SF \enspace}
		\end{array}
		\left\{\begin{array}{c}
			\text{Torsion-free} \\
			\text{classes} 
		\end{array} \right\}
		\]
		and either class determines the other.
		\ii $\ST$ is closed under quotients and extensions, $\SF$ is closed under subobjects and extensions.
		\ii $t(-)$ and $(-)/t(-)$ are functors sending each module $X$ to the corresponding module $tX$ (resp. $X/tX$) in the canonical short exact sequence associated to $X$. The action on morphisms is respectively defined by restriction and factoring. We obtain the following adjunctions:
		\[
		\begin{tikzcd}[column sep=large]
			\ST \rar[bend left, "i"{name=it}, pos=0.55] & A\Mod \lar[bend left, "t(-)"{name=t}] 
			\rar[bend left, "(-)/t(-)"{name=f}, pos=0.4] & \SF \lar[bend left, "i"{name=it2}]
			\arrow[phantom, from=it, to=t, "\dashv" rotate=-90]
			\arrow[phantom, from=f, to=it2, "\dashv" rotate=-90]
		\end{tikzcd}
		\]
	\end{enumerate}
\begin{proof}
	In each statement of the proposition we prove just one of the two results, the other is dual.
\begin{enumerate}[left=0pt]
	\ii We prove that $\SF = \ST^{\perp 0}$. $\SF \subset \ST^{\perp 0}$ is clear. Now take $X \in \ST^{\perp 0}$ and consider its associated short exact sequence \[
	\begin{tikzcd}
		0 \rar & tX \rar["= 0"] & X \rar & X/tX \rar & 0.
	\end{tikzcd}\]
    That $X \in \ST^{\perp 0}$ immediately gives us that $tX = 0$, thus $X \cong X/tX \in \SF$. That these bijections are inclusion-reversing is easy to see: Let $\ST' \subset \ST$. Show that $\ST^{\perp 0} \subset \ST'^{\perp 0}$. Taking any $F \in \ST^{\perp 0}$, clearly $\Hom_A(T, F) = 0$ for all $T \in \ST$. But $\ST' \subset \ST$, so we are done.
	\ii Take a quotient $Y$ of $T \in \ST$ and apply $\Hom_A(-, F)$. \[
	\begin{tikzcd}[column sep=small, row sep=-0.5em]
		T \rar & Y \rar & 0 & \rightsquigarrow & 0 \rar & \Hom_A(Y,F) \rar & \Hom_A(T, F) = 0
	\end{tikzcd}\]
	Allowing $F$ to range across the entirety of $\SF$ we see that $\Hom_A(Y, \SF) = 0$. We now repeat the same procedure for an extension. Consider:
	\[
	\begin{tikzcd}
		0 \rar & T \rar & Y \rar & T' \rar & 0 & \text{ where } T, T' \in \ST
	\end{tikzcd}
	\]
	Applying $\Hom_A(-, F)$ where $F$ ranges across $\SF$ we obtain \[
	\begin{tikzcd}[row sep=-0.5em]
		& = 0 && = 0 \\
		0 \rar & \Hom_A(T',F) \rar & \Hom_A(Y,F) \rar & \Hom_A(T,F)
	\end{tikzcd}
	\] and see that yet again $Y \in {}^{\perp 0}\SF$.
	\ii $t(-)$ acts on morphisms by restriction to $t(X)$, and by the above remark any morphism $t(X) \rightarrow Y$ factors through $t(Y)$. Consequently, $i$ and $t(-)$ are functors. We show that $i(-)$ is left adjoint to $t(-)$ by defining unit and counit natural transformations \[
	\begin{tikzcd}[row sep=-0.5em, column sep=1.8em]
		& \id_{\ST} \ar[r, Rightarrow] & t \circ i &&& & i \circ t \ar[r, Rightarrow] & \id_{A\Mod} \\
		\eta : &&&&& \epsilon : && \\
		& T \ar[r, "\sim"] & T &&& & tX \ar[r, hookrightarrow] & X
	\end{tikzcd}\] and checking the corresponding 'triangle equations', i.e. 
    \[
	\begin{tikzcd}
		\id_i = \epsilon_{i(-)} \circ i(\eta_{-}) &,& \id_{t} = t(\epsilon_{-}) \circ \eta_{t(-)}.
	\end{tikzcd}\]
	That $\eta$ and $\epsilon$ are natural transformations is clear. We check the triangle equations pointwise. For $T \in \ST$ we have $tT = T$, so both $\eta_T : T \rightarrow tT$ and $\epsilon_{iT} : tT \hookrightarrow T$ are the identity of $T$, and thus so is their composition. For $X \in A\Mod$ we have $ttX = tX$, so $\eta_{tX} : tX \rightarrow ttX$ is the identity of $tX$, and so is $t(\epsilon_X)$, the restriction of the inclusion $tX \hookrightarrow X$ to $ttX = tX$. Therefore the triangle equations are fulfilled, and this shows that $i(-)$ is left adjoint to $t(-)$.
	
	If we wish to interpret this adjunction in terms of \[
	\Hom_A(iT, X) \cong \Hom_\ST(T, tX),
	\]
	then the bijection between these Hom-sets merely allows us to alternate between viewing a map $f$ out of $T \in \ST$ as mapping to $X$ or $tX$ respectively. 
\end{enumerate}
\end{proof}
\end{proposition}

\begin{example}
\begin{enumerate}
    \ii[]
	\ii There exist the trivial torsion theories $(0, A\Mod) ,~(A\Mod, 0)$.
	\ii The torsion and torsion-free modules give a torsion theory in the category of $\ZZ$-modules.
    \ii Let $A$ be an algebra of finite representation type. Partition its indecomposable modules into two sets $T,F$ fulfilling $\Hom_A(T,F) = 0$. Then $(\add T,~\add F)$ is a torsion theory in $A\Mod$: If an indecomposable module $X$ fulfills $\Hom_A(X, F) = 0$, then $X \notin F$, since $\Hom_A(X,X) \neq 0$. Thus $X \in T$, which shows that $\add T = {}^{\perp 0}(\add F)$, and dually $\add F = (\add T)^{\perp 0}$. That this suffices will follow from \cref{2TorsionClassEquivalences}.
    \ii Given some representation-infinite algebra one should expect to discover infinitely many torsion classes—in fact, for large classes of algebras the property of being representation-infinite is equivalent to the existence of infinitely many torsion classes. There are exceptions to this rule, however. At the time of writing this is a topic of active research. We provide a simple example—infinitely many of them, in fact—in \cref{7WildTauFiniteAlgebrasExample}.
\end{enumerate}
\end{example}

\begin{example}[Torsion classes for $K(1 \rightarrow 2)$]
We consider the quiver $Q : 1 \rightarrow 2$. We can construct its AR-quiver as follows:
\[
	Q :
	\begin{tikzcd}
		1 \ar[r] & 2
	\end{tikzcd}
	\qquad\qquad \Gamma(Q) :
	\begin{tikzcd}
		\ul{01} \ar[dr] \art && \ol{10} \\
		& \ul{\ol{11}} \ar[ur] & 
	\end{tikzcd}
\]
The following diagram contains a complete list of torsion theories for $KQ$, here depicted as a poset with the arrows denoting \emph{reverse} inclusion. This is the so-called Hasse quiver of the torsion classes of $KQ$.
\[
\begin{array}{c|c}
	\ST & \SF \\ \hline
	01 \oplus 11 \oplus 10 & 0 \\
	11 \oplus 10 & 01 \\
	10 & 11 \oplus 01 \\
	01 & 10 \\
	0 & 01 \oplus 11 \oplus 10 
\end{array} \qquad \quad
\begin{tikzcd}[column sep=-0.2em, row sep=small]
	& 01\oplus11\oplus10 \dlar \ar[ddr] &\\
	11\oplus 10 \dar \ && \\
	10 \drar && 01 \dlar \\
	& 0  &		
\end{tikzcd}
\]
\end{example}

\begin{remark}\label{2TorsionFunctorsNotExact}
    \begin{enumerate}
        \ii[]
        \ii The functors $t(-),~(-)/t(-)$ are in general not exact, and—viewed as endofunctors of $A\Mod$—not even left- (resp. right-)exact. The adjunctions only guarantee that $t(-)$ preserves limits computed in $\ST$, and that $(-)/t(-)$ preserves colimits computed in $\SF$. These need not agree with the kernels and cokernels in $A\Mod$. In fact, $t(-)$ is left-exact if and only if $\ST$ is closed under submodules. Consider the torsion theories $(\add\{11, 10\},~\add\{01\})$ and $(\add\{10\},~\add\{01, 11\})$ of the above example.
        Applying $t(-)$ for the first, and $(-)/t(-)$ for the second torsion theory to the AR-sequence of $KQ$, we obtain
        \[
            \begin{tikzcd}[column sep=small, row sep=small]
                0 \rar & 01 \rar & 11 \rar & 10 \rar & 0 & \overset{t(-)}{\rightsquigarrow} & 0 \rar & 0 \rar & 11 \rar & 10 \rar & 0 \\
                0 \rar & 01 \rar & 11 \rar & 10 \rar & 0 & \overset{(-)/t(-)}{\rightsquigarrow} & 0 \rar & 01 \rar & 11 \rar & 0 \rar & 0
            \end{tikzcd}
        \]
        Neither sequence on the right is exact at $11$. So $t(-)$ is not left-exact, and $(-)/t(-)$ is not right-exact.
        \ii Over a representation-finite quiver any module class $\SC$ is of the form $\add M$, where $M$ is a finite dimensional module containing a summand corresponding to each isomorphism class of the indecomposable modules contained in $\SC$.
    \end{enumerate}
\end{remark}

\begin{proposition}\label{2TorsionClassEquivalences}
	For a module class $\ST$ in $A\Mod$ the following are equivalent:
	\begin{enumerate}
		\ii $\ST$ is a torsion class for some torsion theory in $A\Mod$.
		\ii $\ST = {}^{\perp 0}(\ST^{\perp 0})$
		\ii $\ST = {}^{\perp 0}\SC$ for some module class $\SC$.
		\ii $\ST$ is closed under quotients and extensions.
	\end{enumerate}
	\begin{proof}
        Using that $\ST = {}^{\perp 0} \SF$ from \cref{2Adjoints} makes proving $(1) \Rightarrow (2) \Rightarrow (3)$ straightforward. A closer look at the proof of statement $(2)$ of \cref{2Adjoints} suffices to show that $(3) \Rightarrow (4)$.
		
		It remains to show $(4) \Rightarrow (1)$: We define $\SF = \ST^{\perp 0}$. Given any $X \in A\Mod$, let's construct an exact sequence. Since $X$ is finite-dimensional we can find a submodule $tX \leq X$ so that $tX$ is of maximal dimension and $tX \in \ST$. Now prove that $X/tX \in \SF$. Consider a map $T \rightarrow X/tX$ with image $S/tX$, where $T \in \ST$. Evidently $S/tX \in \ST$. We obtain a short exact sequence
        \[
            \begin{tikzcd}
                0 \rar & tX \rar & S \rar & S/tX \rar & 0.
            \end{tikzcd}
        \]
        Since $tX,~S/tX$ are in $\ST$, so is $S$. This makes $S$ a submodule of $X$ which is also in $\ST$. $tX$ is a maximal such module, so $S = tX$, and we are done. 
	\end{proof}
\end{proposition}

\subsection{Covers and Envelopes}
\begin{definition}
    Consider a map $\theta : X \rightarrow Y$ of $A$-modules. 
    \begin{enumerate}
        \ii We say that $\theta$ is \emph{codomain minimal} if, given some $\alpha \in \End_A(Y)$ such that $\alpha \theta = \theta$, it follows that $\alpha$ necessarily is invertible.
        \ii We say that $\theta$ is \emph{domain minimal} if, given some $\beta \in \End_A(X)$ such that $\theta \beta = \theta$, it follows that $\beta$ is necessarily invertible.
    \end{enumerate}
    In the literature these are commonly referred to as "left minimal" and "right minimal". These terms are confusing, which is why we chose this nomenclature.
\end{definition}

\begin{lemma}\label{2minimalMapsAlsoMinimalOnSummands}
    We obtain the following immediate consequence:
    \begin{enumerate}
        \ii Given some domain minimal map $\theta : X \rightarrow Y$, the restriction of $\theta$ to each summand of $X$ is domain minimal. It follows that $\theta$ is non-zero on each non-zero summand of $X$.

        \ii Given some codomain minimal map $\theta : X \rightarrow Y$, the projection of $\theta$ to each summand of $Y$ is codomain minimal. It follows that the composition of $\theta$ with the projection onto any non-zero summand of $Y$ is non-zero.
    \end{enumerate}
    Note that the second statement does \emph{not} say that $\im \theta$ has non-trivial intersection with each summand of $Y$. For example, for the path algebra of $2 \rightarrow 1 \leftarrow 3$ the diagonal map $S_1 \rightarrow P_2 \oplus P_3$ is codomain minimal, but its image intersects both $P_2$ and $P_3$ trivially.
    \begin{proof}
        We prove the first result, the second is dual. We take some domain minimal map $\theta : X \rightarrow Y$, and a summand $X_0$ of $X$. Use $X_1$ to refer to the complement of $X_0$ in $X$.

        Take an endomorphism $\alpha_0$ of $X_0$ which fulfills that $\restr{\theta}{X_0} \alpha_0 = \restr{\theta}{X_0}$. Then $\alpha :=
        \begin{pmatrix} \alpha_0 & 0 \\ 0 & \id_{X_1} \end{pmatrix}$
         is an endomorphism of $X$ such that $\theta \alpha = \theta$. Thus $\alpha$ is invertible, and so is $\alpha_0$. This shows that $\restr{\theta}{X_0}$ is domain minimal.

        If $\restr{\theta}{X_0} = 0$, then $\alpha_0 = 0$ fulfills $\restr{\theta}{X_0} \alpha_0 = \restr{\theta}{X_0}$, so the zero endomorphism of $X_0$ is invertible, i.e. $X_0 = 0$.
    \end{proof}
\end{lemma}

\begin{lemma}\label{2existenceMinimalMaps}
    Given a map $\theta : X \rightarrow Y$ of (finite-dimensional) $A$-modules
        \begin{enumerate}
            \ii There is a decomposition $Y = Y_0 \oplus Y_1$ such that $\im \theta \subseteq Y_1$, and $\theta : X \rightarrow Y_1$ is codomain minimal.
            \ii There is a decomposition $X = X_0 \oplus X_1$ such that $X_0 \subseteq \ker \theta$, and $\theta : X_1 \rightarrow Y$ is domain minimal.
        \end{enumerate}
    As a consequence, any non-zero map with indecomposable domain (resp. codomain) is already domain (resp. codomain) minimal.
\begin{proof}
We prove the second statement, the first is dual. Given a map $\theta$, take a decomposition $X = X_0 \oplus X_1$ such that $X_0 \subseteq \ker \theta$ is of maximal dimension. Now take any endomorphism $\alpha$ of $X_1$ such that $\theta = \theta\alpha$.

Fitting's Lemma states that there exists a decomposition $X_1 = X_{1,0} \oplus X_{1,1}$, such that $\restr{\alpha}{X_{1,0}}$ is nilpotent, and $\restr{\alpha}{X_{1,1}}$ is invertible. Using that $\theta = \theta\alpha = \theta\alpha^n$ for arbitrarily high $n$ we see that $\theta$ equals zero on $X_{1,0}$. But if this is the case then $\theta$ is zero on $X_0 \oplus X_{1,0}$ and—by maximality of $X_0$—we see that $X_{1,0}$ is the zero module.
\end{proof}
\end{lemma}

\begin{lemma}\label{2sumMinimalMaps}
    Let there be (co)domain minimal $A$-module maps $\theta_k : X_k \rightarrow Y_k$, where $k$ ranges across $1, \dots, n$.

    Then their direct sum $\theta : \bigoplus_{k=1}^n X_k \rightarrow \bigoplus_{k=1}^n Y_k$ is (co)domain minimal.
\begin{proof}
    We prove this for domain minimality, the other proof is dual. Write $X := \bigoplus_k X_k$, and $\pi_k : X \rightarrow X_k$, $\iota_k : X_k \rightarrow X$ for the canonical projections and inclusions. Since $\theta$ is a direct sum of maps, we have $\ker \theta = \bigoplus_k \ker \theta_k$.

    If $\theta$ is not domain minimal then, by \cref{2existenceMinimalMaps}, there exists a non-zero summand of $X$ contained in $\ker \theta$, and thus also an indecomposable such summand $X'$. Write $\iota : X' \rightarrow X$ and $\pi : X \rightarrow X'$ for the inclusion and projection. Then $\id_{X'} = \sum_k (\pi \iota_k)(\pi_k \iota)$, and since $\End_A(X')$ is local, one of these summands is invertible. For this value of $k$, the map $\pi_k \iota : X' \rightarrow X_k$ is a split monomorphism, so its image is a non-zero summand of $X_k$. It is contained in $\pi_k(\ker \theta) = \ker \theta_k$, so $\theta_k$ is zero on a non-zero summand of $X_k$. This contradicts \cref{2minimalMapsAlsoMinimalOnSummands}, and thus $\theta$ is domain minimal.
\end{proof}
    We note that a converse of this statement also holds, but we don't need it for our purposes. We refer the reader to \cite{NA2}, which proves both directions using a different approach.
\end{lemma}

\begin{definition}[Covers and Envelopes]
Take a module class $\SC$ and some module $X$ in $A\Mod$. 
\begin{enumerate}
	\ii A $\SC$-\emph{preenvelope} of $X$ is a map $\theta : X \rightarrow C$, where $C \in \SC$, so that every morphism from $X$ to a module $C'$ in $\SC$ factors through $\theta$.
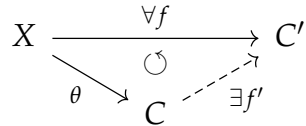
	\[
        \begin{tikzcd}[row sep=small]
		X \ar[rr, "\forall f"{name=U}] \ar[dr, "\theta", swap] && C'\\
		& C \ar[ur, dashed, "\exists f'", swap] \ar["\circlearrowleft", from=U, phantom] & {}
	\end{tikzcd}
	\]
	\ii A \emph{$\SC$-envelope} is a \emph{codomain minimal} $\SC$-preenvelope.
	
	\ii A $\SC$-\emph{precover} of $X$ is a map $\theta : C \rightarrow X$, where $C \in \SC$, so that every morphism from a module $C'$ in $\SC$ to $X$ factors through $\theta$.
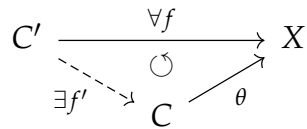
	\[
        \begin{tikzcd}[row sep=small]
		C' \ar[rr, "\forall f"{name=U}] \ar[dr, dashed, "\exists f'", swap] && X \\
		& C \ar[ur, "\theta", swap] \ar["\circlearrowleft", from=U, phantom] & {}
	\end{tikzcd}
    \]
	\ii A \emph{$\SC$-cover} is a \emph{domain minimal} $\SC$-precover.
\end{enumerate}
    In the literature $\SC$-(pre)envelopes are sometimes referred to as "left $\SC$-approximations". Similarly, $\SC$-(pre)covers are sometimes referred to as "right $\SC$-approximations".
\end{definition}

\begin{proposition}\label{2coversProperties}
\begin{enumerate}
    \ii[]
	\ii If $X$ has a $\SC$-preenvelope then it also has a $\SC$-envelope.
    \ii[1')] If $X$ has a $\SC$-precover then it also has a $\SC$-cover.
	\ii Any two $\SC$-envelopes of $X$ are isomorphic.
    \ii[2')] Any two $\SC$-covers of $X$ are isomorphic.
    \ii Given $n$ $\SC$-envelopes $X_k \xrightarrow{\quad \theta_k\quad} C_k$, their direct sum $\bigoplus_k X_k \xrightarrow{\bigoplus_k \theta_k} \bigoplus_k C_k$ is also an $\SC$-envelope.
    \ii[3')] Given $n$ $\SC$-covers $C_k \xrightarrow{\quad \theta_k\quad} X_k$, their direct sum $\bigoplus_k C_k \xrightarrow{\bigoplus_k \theta_k} \bigoplus_k X_k$ is also an $\SC$-cover.
\end{enumerate}
\begin{proof}
    Naturally, these two cases are dual to each other.
\begin{enumerate}
    \ii This follows from \cref{2existenceMinimalMaps}.
	\ii Consider $\SC$-envelopes $\theta : X \rightarrow C$, $\theta' : X \rightarrow C'$. $\theta,~\theta'$ respectively induce maps $f : C \rightarrow C'$, $g : C' \rightarrow C$, such that $f\theta = \theta'$ and $g\theta' = \theta$. Now $gf \in \End_A(C)$ and $fg \in \End_A(C')$ fulfill $gf\theta = \theta$ and $fg\theta' = \theta'$. By codomain minimality $fg$,~$gf$ are both invertible, and thus $f$ and $g$ are isomorphisms.
    \ii A direct sum of $\SC$-preenvelopes is an $\SC$-preenvelope, since $\SC$ is closed under direct sums. The claim now follows from \cref{2sumMinimalMaps}.
\end{enumerate}
\end{proof}
\end{proposition}

\begin{definition}
A module class $\SC$ in $A\Mod$ is:
\begin{enumerate}
	\ii \emph{covariantly finite} if every $A$-module has a $\SC$-envelope.
	\ii \emph{contravariantly finite} if every $A$-module has a $\SC$-cover.
	\ii \emph{functorially finite} if it is both covariantly finite and contravariantly finite.
\end{enumerate}
\end{definition}

\begin{example}
    For any two modules $M, X$: Taking a spanning set of $\Hom_A(X,M)$ and considering the corresponding map $X \rightarrow M^n$ provides us with an $\add M$-preenvelope of $X$.
	
	Similarly, taking a spanning set of $\Hom_A(M,X)$ and considering the corresponding map $M^n \rightarrow X$ provides us with an $\add M$-precover of $X$.

    By \cref{2coversProperties} the existence of preenvelopes and precovers implies the existence of envelopes and covers, and therefore $\add M$ is functorially finite, as long as $M$ is a finite dimensional module.

    From this it immediately follows that over a representation-finite algebra every module class is functorially finite.
\end{example}

\begin{example}
Take some module class $\SC$, and consider the inclusion functor $i : \SC \hookrightarrow A\Mod$. If there exists a left adjoint functor $L \dashv i$, then the natural unit map $\eta : \id_{A\Mod} \Rightarrow iL$ is a $\SC$-envelope at each component.

If there exists a right adjoint functor $i \dashv R$, then the natural counit map $\epsilon : iR \Rightarrow \id_{A\Mod}$ is a $\SC$-cover at each component.

For example:
\begin{enumerate}
    \ii Given a torsion theory $(\ST, \SF)$, $\ST$ is contravariantly finite, and $\SF$ is covariantly finite. Given any module $X$ in $A\Mod$, $ t(X) \hookrightarrow X $ is a $\ST$-cover of $X$, and $ X \twoheadrightarrow X/t(X) $ is a $\SF$-envelope of $X$.
    \ii $X \rightarrow \text{top}(X) := X/\rad X$ is a semisimple envelope, and $\text{socle}(X) \rightarrow X$ is a semisimple cover.
\end{enumerate}
A warning: Projective covers are $\add {}_AA$-covers, injective envelopes are $\add D(A_A)$-envelopes, but neither is an instance of the above pattern.

\begin{proof}
Take the path algebra $K(1 \rightarrow 2)$, and denote the injective envelope with $E(-)$, and the inclusion of the (functorially finite) module class of injective modules into $KQ\Mod$ with $i(-)$. Assume that $E \dashv i$. In this case the maps $X \hookrightarrow E(X)$ assemble into a natural transformation. Now consider the following diagram, obtained by applying $E(-)$ to the AR-sequence in $KQ\Mod$.
\[
    \begin{tikzcd}[column sep=large]
    0 \rar & 01 \rar["f"] \dar[hook] & 11 \rar["g"] \dar[hook] & 10 \rar \dar[hook] & 0 \\
           & \underset{=11}{E(01)} \rar["Ef"] & \underset{=11}{E(11)} \rar["Eg"] & \underset{=10}{E(10)} & 
\end{tikzcd}
\]
The commutativity of the first square implies that $Ef = \id_{11}$, and the commutativity of the second square implies that $Eg \neq 0$. This is a contradiction since the commutativity of the outer rectangle requires their composition to be the zero morphism.
\end{proof}
\end{example}

We have seen that every torsion class $\ST$ which can be written in the form $\ST = \add T$ for some module $T$ is functorially finite. This is a very restrictive condition. In the following we will look at torsion classes which are \emph{generated} by a module $T$, i.e. can be written as a quotient of a module in $\add T$.

We will see that any torsion class that can be written in the form $\gen T$ is already functorially finite. In the next section we will prove a converse to this result: Every functorially finite torsion class can be written in the form $\gen T$, where $T$ is a module with very interesting properties.

\begin{definition}
    Given some module $M$ in $A\Mod$ we write $\gen M$ for the class of modules generated by $M$, i.e. those modules $X$ which can be written as a quotient of a module $M_0$ in $\add M$:
    \[
        \begin{tikzcd}
            M_0 \rar & X \rar & 0
        \end{tikzcd}
    \]

    We write $\gen_1 M$ for the class of modules with an $\add M$-presentation, i.e. those modules $X$ for which there exists an exact sequence
    \[
        \begin{tikzcd}
            M_1 \rar & M_0 \rar & X \rar & 0
        \end{tikzcd}
    \]
    with $M_1, M_0$ in $\add M$. These definitions dualise as $\cogen M$, and $\cogen_1 M$.
\end{definition}

\begin{lemma}\label{2genMCovariantlyFinite}
For any module $M$: $\gen M$ is covariantly finite. Consequently, if $\gen M$ is a torsion class, then it is functorially finite.
\begin{proof}
Take any module $X$, and a projective cover $P$ of $X$. Next take an $\add M$-preenvelope $M^n$ of $P$, and form the pushout: \[
\begin{tikzcd}
P \rar[two heads] \dar \drar[phantom, "\ulcorner", near start] & X \dar \\
M^n \rar[two heads] & G
\end{tikzcd}
\]
The surjectivity of $M^n \twoheadrightarrow G$ follows from the surjectivity of $P \twoheadrightarrow X$. We take some map $f : X \rightarrow G'$, where $G' \in \gen M$, and prove that there exists a map $M^n \rightarrow G'$ which agrees with $f$ on $P$. There is a map $M^m \twoheadrightarrow G'$, and the projectivity of $P$ ensures that $P \twoheadrightarrow X \rightarrow G'$ factors through $M^m$. Now $M^n$ is an $\add M$-preenvelope of $P$, so we obtain an induced map $M^n \rightarrow M^m \twoheadrightarrow G'$.
\[
\begin{tikzcd}
	P \rar[two heads] \dar \drar[phantom, "\ulcorner", near start] & X \dar \ar[ddr, bend left, "f"] & \\
	M^n \drar \rar[two heads] & G \drar[dashed, "\exists"] & \\
	& M^m \rar[two heads] & G'
\end{tikzcd}
\]
Since these maps agree on $P$ we obtain an induced map $G \rightarrow G'$, and thus $X \rightarrow G$ is a $\gen M$-preenvelope of $X$. By \cref{2coversProperties}, $X$ also has a $\gen M$-envelope.
\end{proof}
\end{lemma}

Given some module $M$, we would like to know under which circumstances $\gen M$ is a torsion class in the first place. For reasons of conceptual clarity we prepend the following lemma:

\begin{lemma}\label{2genMWhenTorsionClass}
    If $ \Ext^1_A(M, \gen M) = 0 $ holds, then $\gen M$ is a torsion class.
    \begin{proof}
     Evidently $\gen M$ is closed under quotients, and by \cref{2TorsionClassEquivalences} we merely need to figure out when $\gen M$ is closed under extensions. So, let us take an extension of $X, Y \in \gen M$, and then take the natural pullback of $M^n \twoheadrightarrow Y$.  
	\[
	\begin{tikzcd}
		0 \rar & X \rar \dar[equal] & E' \rar \dar \drar["\lrcorner", phantom] & M^n \rar \dar[two heads] & 0 \\
		0 \rar & X \rar & E \rar & Y \rar & 0
	\end{tikzcd}
	\]
	$E' \rightarrow E$ is surjective, therefore $E' \in \gen M$ already shows that $E \in \gen M$. This reduces the difficulty of the situation by allowing us to restrict to extensions of the form $\Ext_A^1(M, X)$ for $X \in \gen M$, and it suffices if all of these split. This is what we wanted to show.
\end{proof}
\end{lemma}

A major result known as the \emph{Auslander-Smal\o\ Lemma}, which we will prove momentarily, shows that
\[
    \Ext^1_A(M, \gen M) = 0 \qquad \iff \qquad \Hom_A(M, \tau M) = 0.
\]
This statement can be viewed as the starting point of \ta-tilting theory. Modules which fulfill these equivalent conditions are referred to as \emph{\ta-rigid modules}.

We have just shown that any \ta-rigid module generates a functorially finite torsion class $\gen M$. At the end of the next section we will prove a converse of this statement: Any functorially finite torsion class arises from a \ta-rigid module.

While that is interesting in its own right, we will see throughout our journey that \ta-rigid modules possess many more interesting properties, and assemble into rich combinatorial structures.

\section{\texorpdfstring{$\tau$}{Tau}-rigid Modules}
\subsection{The Auslander-Smal\o\ Lemma and Theorem}

\begin{definition}
A module $M$ in $A\Mod$ is called \emph{\ta-rigid} if \[
\Hom_A(M, \tau M) = 0. \]
\end{definition}

Our first goal is to prove the aforementioned Auslander-Smal\o\ Lemma. We will then explore its results, and prove that every functorially finite torsion class arises from a \ta-rigid module.

\begin{lemma}[Auslander-Smal\o]\label{3AuslanderSmaloLemma}
	Let $M, N$ be $A$-modules. The following are equivalent:
	\begin{enumerate}
		\ii $\Hom_A(N, \tau M) = 0$
		\ii $\Ext_A^{1}(M, \gen N) = 0$
	\end{enumerate}
\begin{proof}
$1)\implies2):$ Take a module $N'$ in $\gen N$. There exists a surjection $N^n \rightarrow N' \rightarrow 0$, and we apply $\Hom_A(-, \tau M)$ to obtain
\[
    \begin{tikzcd}
        0 \rar & \Hom_A(N', \tau M) \rar & \Hom_A(N^n, \tau M) = 0
    \end{tikzcd}.
\]
Using the Auslander-Reiten-Formula it follows that
\[
    D\Ext^1_A(M, N') \cong \ol{\Hom}_A(N', \tau M) = 0 \text{\qquad for all } N' \text{ in } \gen N ,
\]
and this is up to duality exactly what we needed to show.

$2)\implies 1):$ Prove the contrapositive. Take a non-zero map $f : N \rightarrow \tau M$ and factorize it as follows:
\[
\begin{tikzcd}
	N \rar[twoheadrightarrow] & G \rar["h", hook] & \tau M
\end{tikzcd}
\]
Clearly $G$ is in $\gen N$. It suffices to show that $\Ext_A^1(M,G) \neq 0$, this is (again by the AR-formula) equivalent to showing that $\ol{\Hom}(G, \tau M) \neq 0$, and for this it is sufficient to prove that $h$ cannot factor through an injective module. 

Therefore, assume that $h$ factors through an injective module. In this case it will factor through the injective envelope $E(G)$ of $G$.
\[
\begin{tikzcd}
	G \ar[rr, "h", bend left] \rar[hook] & E(G) \rar["k"] & \tau M
\end{tikzcd}\]
If $k$ were injective it'd be a split-mono, and since $\tau M$ has no injective summand this cannot be the case. $G$ is essential in $E(G)$, i.e. meets every submodule, including $\ker k$. It follows that $\ker h \neq 0$. This is a contradiction to the injectivity of $h$, therefore $h$ cannot factor through an injective module. 
\end{proof}
\end{lemma}

\begin{definition}
Given a module class $\SC$ in $A\Mod$ and a module $M$ contained in $\SC$, we say that 
	\begin{enumerate}
		\item $M$ is \emph{Ext-projective} in $\SC$ if $\Ext_A^{1}(M,\SC) = 0$.
		\item $M$ is \emph{Ext-injective} in $\SC$ if $\Ext_A^{1}(\SC, M) = 0$.
	\end{enumerate}
Naturally, the Ext-projectives of $A\Mod$ are the projective modules. A complete classification of the Ext-injectives of $A\Mod$ is left as an exercise for the reader.
\end{definition}

\begin{lemma}\label{3extProjectivesBijection}
	If $(\ST, \SF)$ is a torsion theory in $A\Mod$, then 
	\begin{enumerate}
		\item $X \in \ST$ is Ext-projective for $\ST$ iff $\tau X \in \SF$.
		\item $X \in \SF$ is Ext-injective for $\SF$ iff $\tau^{-} X \in \ST$.
		\ii There are bijections
		\[ \left\{
		\begin{array}{c}
			\text{Non-projective} \\
			\text{indecomposable} \\
			\text{Ext-projectives in $\ST$}  \\
			\text{up to isomorphism.}
		\end{array} \right\}
		\begin{array}{c}
			\xrightarrow{\enspace \, \tau \,\enspace} \\
			\xleftarrow{\enspace \tau^- \enspace}
		\end{array}
		\left\{\begin{array}{c}
			\text{Non-injective} \\
			\text{indecomposable} \\
			\text{Ext-injectives in $\SF$} \\
			\text{up to isomorphism.}
		\end{array} \right\}
		\]
\end{enumerate}
\begin{proof}
\begin{enumerate}[left=0pt]
    \ii[]
\ii Using the Auslander-Smal\o{} Lemma, Ext-projectivity of $X \in \ST$ for $\ST$ is equivalent to $\Hom_A(T, \tau X) = 0$ for all $T \in \ST$, i.e. $\tau X \in \ST^{\perp 0} = \SF$.
\ii Using the dual of the Auslander-Smal\o\ lemma, Ext-injectivity of $X \in \SF$ is equivalent to $\Hom_A(\tau^- X, F) = 0$ for all $F \in \SF$, i.e. $\tau^- X \in {}^{\perp 0}\SF = \ST$. The dual of the Auslander-Smal\o\ lemma is stated in \cref{5dualAuslanderSmaloLemma}, in a subsection in which we take a closer look at the theory dual to the one which we are developing.
\ii Follows: Just remember that $\tau, \tau^-$ give bijections between non-projective indecomposables up to isomorphism, and non-injective indecomposables up to isomorphism.
\end{enumerate}
\end{proof}
\end{lemma}

\begin{lemma}
Given a torsion theory $(\ST, \SF)$ in $A\Mod$ we have that
\begin{enumerate}
	\ii The Ext-injectives for $\ST$ are the modules $tI$ with $I$ injective. 

    If $I$ is an indecomposable injective module which is not contained in $\SF$, then $tI$ is an indecomposable Ext-injective, and every indecomposable Ext-injective for $\ST$ arises this way.

	\ii The Ext-projectives for $\SF$ are the modules $P/tP$ with $P$ projective. 

    If $P$ is an indecomposable projective module which is not contained in $\ST$, then $P/tP$ is an indecomposable Ext-projective, and every indecomposable Ext-projective for $\SF$ arises this way.
\end{enumerate}
\begin{proof}
We prove the first statement, the second is dual. Consider a short exact sequence
\[
\begin{tikzcd}
0 \rar & tI \rar & E \rar & T \rar & 0
\end{tikzcd}
\]
in $\ST$, where $I$ is some injective module in $A\Mod$. Take the pushout along the canonical inclusion $tI \hookrightarrow I$. The lower short exact sequence splits, since $I$ is injective.
\[
\begin{tikzcd}
0 \rar & tI \rar \dar[hook] \drar[phantom, "\ulcorner", near start] & E \rar \dar & T \rar \dar[equal] & 0 \\
0 \rar & I \rar & E' \rar & T \rar & 0
\end{tikzcd}
\]
Composing $E \rightarrow E'$ with a retraction $E' \rightarrow I$ of the lower sequence, we obtain a map $E \rightarrow I$ which restricts to the inclusion $tI \hookrightarrow I$. Since $E \in \ST$, this map factors through $tI$, and we obtain a retraction $E \rightarrow tI$ of the upper sequence. Thus $tI$ is Ext-injective.

Conversely, suppose that $X$ is an Ext-injective module and take its injective envelope $E(X)$. The map $X \hookrightarrow E(X)$ is in fact a map $X \hookrightarrow tE(X)$. 
\[
\begin{tikzcd}
0 \rar & X \rar & tE(X) \rar & tE(X)/X \rar & 0
\end{tikzcd}
\]
splits, therefore $X$ is a summand of $tE(X)$. Since $X$ is essential in $E(X)$ it follows that $X = tE(X)$. 

It remains to deal with indecomposability. An indecomposable injective module $I$ has a simple socle, so every non-zero submodule of $I$ has a simple socle, and is in particular indecomposable. Thus $tI$ is indecomposable if and only if $tI \neq 0$, which is the case if and only if $I \notin \SF$. Conversely, if $X$ is an indecomposable Ext-injective, then $X = tE(X) = \bigoplus_j tI_j$, where $E(X) = \bigoplus_j I_j$ is a decomposition into indecomposable injectives. Since $X$ is indecomposable, $X = tI_j$ for some $j$.
\end{proof}
\end{lemma}

\begin{example}
	Consider $K(1 \rightarrow 2 \rightarrow 3)$. Define $M = 010 \oplus 111$, which is a \ta-rigid $KQ$-module. Then $\ST := \gen M = \add \{ 010, 111, 110, 100 \}$, and $\SF = \add \{ 001, 011 \}$. 
	
	The Ext-projectives for $\ST$ are $010, 111, 110$. The Ext-injectives for $\ST$ are $111, 110, 100$.
	
	The Ext-projectives for $\SF$ are $001, 011$. The Ext-injectives for $\SF$ are $001, 011$.
	
    The torsion theory can be visualised in the AR-quiver of $KQ$:
    \[
		\Gamma(KQ) \qquad : \qquad 
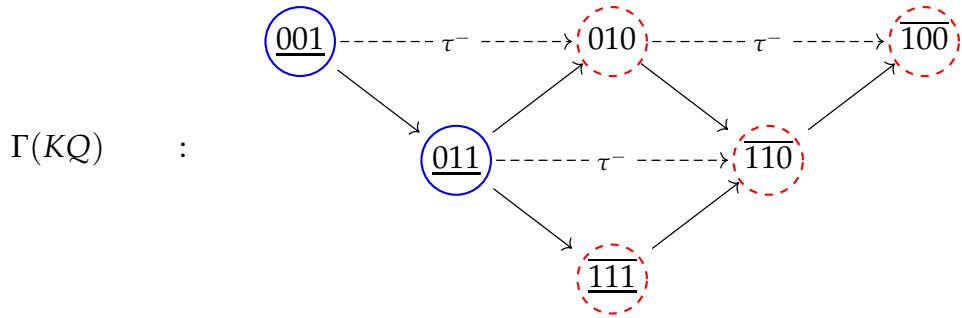
		\begin{tikzcd}[execute at end picture={
				\draw[thick, red, dashed] (\tikzcdmatrixname-2-4) circle[radius=13pt];
				\draw[thick, red, dashed] (\tikzcdmatrixname-1-3) circle[radius=13pt];
				\draw[thick, red, dashed] (\tikzcdmatrixname-3-3) circle[radius=13pt];
                \draw[thick, red, dashed]  (\tikzcdmatrixname-1-5) circle[radius=13pt];
				\draw[thick, blue] (\tikzcdmatrixname-1-1) circle[radius=13pt];
				\draw[thick, blue] (\tikzcdmatrixname-2-2) circle[radius=13pt];
			}]
			\ul{001}\ar[dr] \art && 010 \ar[dr] \art && \ol{100} \\
			& \ul{011} \ar[dr] \ar[ur] \art && \ol{110} \ar[ur] & \\
			&& \ul{\ol{111}} \ar[ur] &&
		\end{tikzcd}	
    \]
\end{example}

\begin{proposition}\label{3tauRigidModulesGenerateTorsionClasses}
Given some $A$-module $T$, the following are equivalent:
\begin{enumerate}
	\ii $T$ is \ta-rigid.
	\ii $\Ext_A^1(T, \gen T) = 0$.
	\ii $\gen T$ is a torsion class and $T$ is Ext-projective in $\gen T$.
	\ii $T$ is Ext-projective in some torsion class.
\end{enumerate}
\begin{proof}
$1) \iff 2)$: Given by the Auslander-Smal\o{} Lemma. 

$2) \implies 3)$: That $\gen T$ is a torsion class is given by \cref{2genMWhenTorsionClass}, and $T$ is Ext-projective in $\gen T$ by $2)$.

$3) \implies 4)$: Trivial.

$4) \implies 2)$: If $T$ is Ext-projective in some torsion class $\ST$, then $\Ext_A^1(T, S) = 0$ for all $S \in \ST$. Since $\gen T$ is a subset of $\ST$, we are done.
\end{proof}
\end{proposition}

\begin{definition}\label{3defExtProjectives}
	Let $(\ST, \SF)$ be a torsion theory over some algebra $A$. We define $\SP(\ST)$ as the direct sum consisting of the indecomposable Ext-projective modules in $\ST$, choosing one summand for each isomorphism class.

    This is a finite direct sum, and thus a well-defined module: Any finite direct sum of pairwise non-isomorphic indecomposable Ext-projectives in $\ST$ is itself Ext-projective in $\ST$, and hence \ta-rigid by \cref{3tauRigidModulesGenerateTorsionClasses}. In \cref{4tauRigidUpperBound} we will see that a \ta-rigid module has at most $\# A$ pairwise non-isomorphic indecomposable summands. Thus $\SP(\ST)$ has at most $\# A$ summands, even if $\ST$ is not functorially finite.
\end{definition}

\begin{example}\label{3LinearA3WithRelationExample}
    Consider the algebra $A$ given by the linearly oriented $\mathbb{A}_3$-quiver with relation. We obtain the following AR-quiver:
	\[
        Q : 
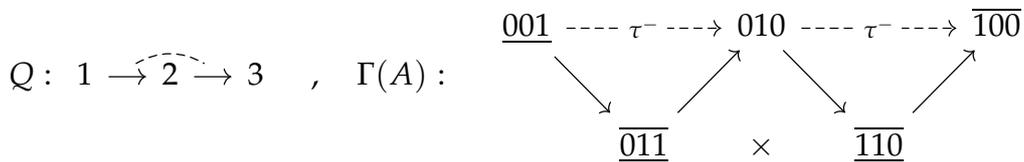
\begin{tikzcd}[column sep=small]
		1 \ar[r, ""{name=a}] & 2 \ar[r, ""{name=b}] & 3
		\arrow[dashed, no head, bend left, from=a, to=b]
    \end{tikzcd} \quad , \quad \Gamma(A) : \quad
	\begin{tikzcd}[column sep=small]
		\ul{001} \art \ar[dr] && 010 \ar[dr] \art && \ol{100}\\
		& \ol{\ul{011}} \ar[rr, phantom, "\times" marking] \ar[ur] && \ol{\ul{110}} \ar[ur] &
	\end{tikzcd}
	\]
	We classify the torsion classes over $A$, and draw the corresponding \emph{Hasse quiver}:
	
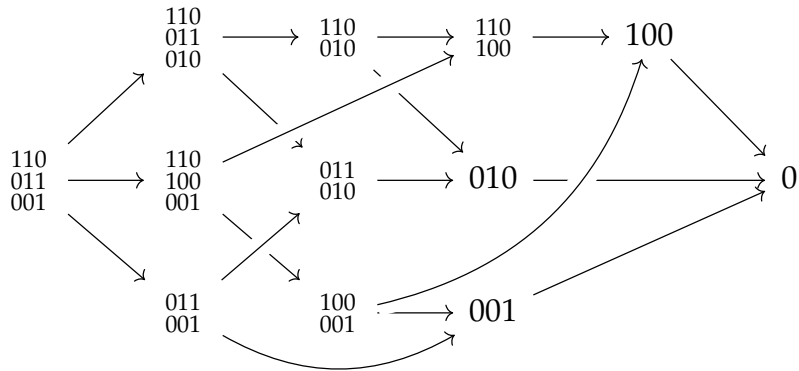
\[
	\begin{tikzcd}
		& \indd{110}{011}{010} \ar[r] \ar[dr] 
		& \ind{110}{010} \ar[r] \ar[dr]
		& \ind{110}{100}  \ar[r]
		& 100 \ar[dr]
		& \\		
		\indd{110}{011}{001} \ar[ur] \ar[r] \ar[dr] 
		& \indd{110}{100}{001} \ar[urr, crossing over] \ar[dr] 
		& \ind{011}{010} \ar[r]
		& 010  \ar[rr]
		&
		& 0 \\
		& \ind{011}{001} \ar[ur, crossing over] \ar[rr, bend right] 
		& \ind{100}{001} \ar[r] \ar[uurr, bend right, crossing over]
		& 001 \ar[urr]
		&
		& \\
	\end{tikzcd}
	\] 
	The vertices of this poset-graph designate (generators of) torsion classes over $A$. Since $A$ is representation-finite, every torsion class $\ST$ can be written as $\ST = \gen \SP(\ST)$.

	The arrows denote reverse inclusion. There are several things the reader might observe here:
	\begin{enumerate}
		\ii Each vertex of the graph has a total of three edges connected to it. Along each edge exactly one module is replaced or removed.
        \ii The support rank\footnote{The support rank of a module class $\SC$ is the number of vertices where the module class is non-zero. Alternatively, the number of isomorphism classes of indecomposable projectives $P$ such that $\Hom_A(P, \SC) \neq 0$.} of each torsion class equals the number of summands of $\SP(\ST)$.
	\end{enumerate}
\end{example}

\begin{example}\label{3ExampleSkewedTriangle}
	For a more complex example, consider the algebra defined as $A := KQ/J$ where
	\[
	\begin{tikzcd}
		Q : & 1 \ar[r, ""{name=a}] \ar[rr, bend right] & 2 \ar[r, ""{name=b}] & 3 \ar[from=a, to=b, bend left, dashed, no head]
	\end{tikzcd}
	\]
	and where $J$ is generated by the composition denoted by the dashed line.
    From the get-go it is clear that this algebra cannot have a preprojective component. This can be seen by noting that there are two non-isomorphic indecomposables with the same dimension vector.\footnote{If it were possible to construct this algebra with a straightforward knitting procedure, then a classic result says that every indecomposable has to be directing, and is already uniquely determined by its dimension vector. See \cite{zbMATH07150963}, Theorem IV.2.12.} The two indecomposables are:
	\[
	\begin{tikzcd}[row sep=tiny]
		& \ul{111} := &&&& \ol{111} := & \\
		K \rar["1"] \ar[rr, bend right, "1"] & K \rar["0"] & K &&
		K \rar["0"] \ar[rr, bend right, "1"] & K \rar["1"] & K
	\end{tikzcd}
	\]
	
	We obtain the AR-quiver through liberal use of dark magic:
	
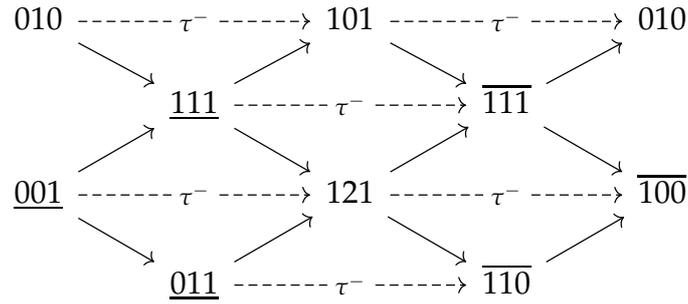
\[
	\begin{tikzcd}[row sep=small]
		010 \art \drar && 101 \drar \art && 010 \\
		& \ul{111} \urar \drar \art && \ol{111} \urar \drar \\
		\ul{001} \urar \drar \art && 121 \urar \drar \art && \ol{100} \\
		& \ul{011} \urar \art && \ol{110} \urar
	\end{tikzcd}
	\]
	Here it should be noted that the "two" indecomposables denoted by $010$ are in fact one and the same. We compute the Hasse quiver of the torsion classes of $A$, and denote a torsion class $\ST$ by the dimension vectors of the summands of $\SP(\ST)$. This is unambiguous since $\ol{111}$ is \emph{not} \ta-rigid, meaning that it cannot possibly occur here.

	
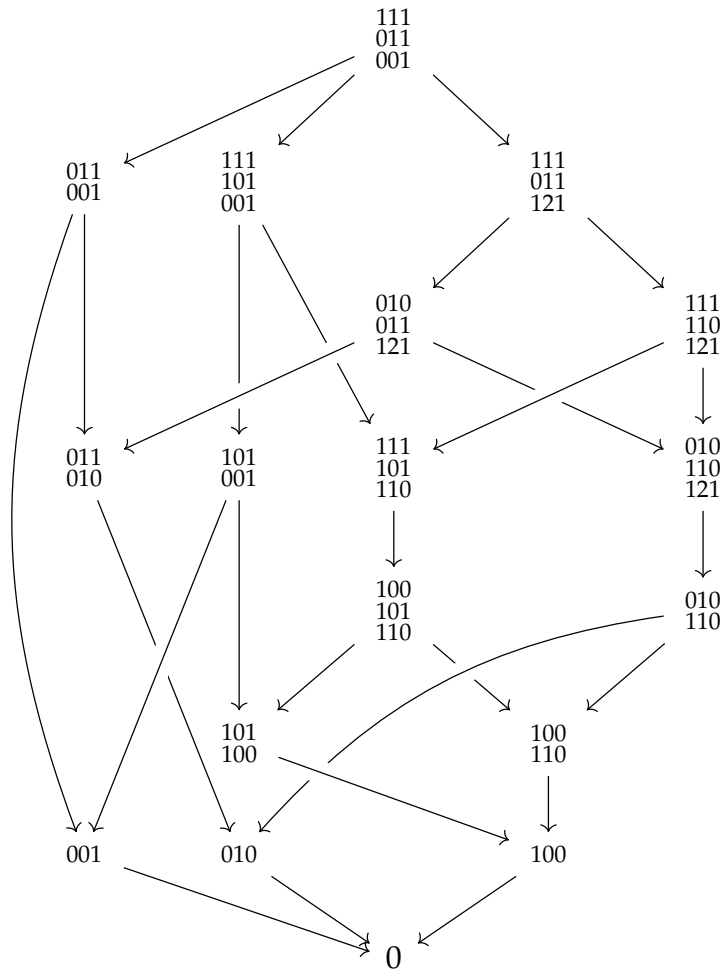
\[
	\begin{tikzcd}[row sep=normal]
		&& \indd{111}{011}{001} \ar[dll] \ar[dl] \ar[dr] &&& \\
		\ind{011}{001} \ar[dd] \ar[ddddd, bend right=20] & \indd{111}{101}{001} \ar[ddr, crossing over] \ar[dd, crossing over] && \indd{111}{011}{121} \ar[dl] \ar[dr] && \\
		&& \indd{010}{011}{121} \ar[dll, crossing over] \ar[drr, crossing over] && \indd{111}{110}{121} \ar[d] \ar[dll, crossing over] & \\
		\ind{011}{010} \ar[dddr, crossing over] & \ind{101}{001} \ar[dddl, crossing over] \ar[dd] & \indd{111}{101}{110} \ar[d] && \indd{010}{110}{121} \ar[d] & \\
		&& \indd{100}{101}{110} \ar[dl] \ar[dr] && \ind{010}{110} \ar[dl] \ar[ddlll, bend right=20, crossing over] & \\
		& \ind{101}{100} \ar[drr, crossing over] && \ind{100}{110} \ar[d] && \\
		\begin{smallmatrix}001\end{smallmatrix} \ar[drr] & \begin{smallmatrix}010\end{smallmatrix} \ar[dr] && \begin{smallmatrix}100\end{smallmatrix} \ar[dl] && \\
		&& 0 &&&
	\end{tikzcd}
	\]
\end{example}

\begin{example}[Kronecker quiver]\label{3KroneckerFFTorsionClasses}
    Consider the Kronecker-quiver $\begin{tikzcd}
    Q: 1 \ar[r, shift right] \ar[r, shift left] & 2 \end{tikzcd}$. We start knitting to obtain its preprojective component $\mathbb{P}$, denoting modules by their dimension vectors:
    \[
        \mathbb{P} : \quad
        \begin{tikzcd}[column sep=small, row sep=small]
         \ul{01}   \art \drar[shift left] \drar[shift right] &&  23   \art \drar[shift left] \drar[shift right]  && \dots\\
        &  \ul{12}  \art \urar[shift left] \urar[shift right] &&  34  \urar[shift left] \urar[shift right] &
    \end{tikzcd}
    \]
    For each indecomposable module $T$ in $\mathbb{P}$ we obtain a torsion class $\gen T$. To see this we use \cref{3tauRigidModulesGenerateTorsionClasses}, and note that any such $T$ is \ta-rigid. Applying $\gen(-)$ to the indecomposables in $\mathbb{P}$ gives us an infinite descending family of functorially finite torsion classes.\footnote{The case when $T$ equals the unique simple projective module is a special case.}
\end{example}

We already have several results concerning the connection between functorially finite torsion classes, \ta-rigid modules, and torsion classes of the form $\gen T$. Let us summarize:

\begin{remark}\label{3tauRigidsGenerateFFTorsionClasses}
    Given some \ta-rigid module $T$, $\gen T$ is a functorially finite torsion class.
    \begin{proof}
        By \cref{2genMCovariantlyFinite} we know that any torsion class, which can be written as $\gen T$ for some module $T$, is functorially finite. By \cref{3tauRigidModulesGenerateTorsionClasses} we know that if we are given some \ta-rigid module $T$, then it is guaranteed that $\gen T$ is in fact closed under extensions, and thus a torsion class.
    \end{proof}
\end{remark}

The following result provides the promised converse: It shows that, given some functorially finite torsion class $\ST$, we can find a \ta-rigid module $T$ such that $\ST = \gen T$. If we can find such a module at all, then $\SP(\ST)$ will always fulfill the same property. This proof is from \cite{NA2}, like most from this section.

\begin{theorem}[Auslander-Smal\o]\label{AuslanderSmaloTheorem}
Let $\ST$ be a functorially finite torsion class. Denote a $\ST$-envelope of $A$ with $f : A \rightarrow T$. Write $T' = T \oplus \coker(f)$. Then:
\begin{enumerate}
	\ii $\ST = \gen T = \gen T'$.
	\ii $T$ is a splitting projective in $\ST$, i.e. for $X \in \ST$ any surjective map $X \twoheadrightarrow T$ splits.
	\ii $T, T'$ are Ext-projective in $\ST$, so they are \ta-rigid.
	\ii Any module $X \in \ST$ is a quotient of a module in $\add T'$ by a submodule in $\ST$.
	\ii Any Ext-projective in $\ST$ is in $\add T'$.
\end{enumerate}
This theorem shows that the Ext-projectives of a functorially finite torsion class possess properties similar to those of projective modules in $A\Mod$. The trivial case $\ST = A\Mod$ makes these similarities even more obvious.
\begin{proof}
\begin{enumerate}
    \ii[]
\ii Given some $X \in \ST$ take a map $A^{n} \twoheadrightarrow X$. Each component factors through $T$, so $T^{n} \twoheadrightarrow X$. 
\ii \[
\begin{tikzcd}
	A \drar["f"] \dar[dashed, "\exists", swap] &\\
	X \rar["\alpha", two heads] & T
\end{tikzcd} \qquad \implies \qquad
\begin{tikzcd}
	X \rar["\alpha", two heads, shift left] & T \lar[dashed, "\exists \beta", shift left] \ar[loop right, "\alpha \beta"]
\end{tikzcd} 
\]
We take some surjection $\alpha : X \twoheadrightarrow T$ as our starting point. Projectivity of $A$ shows that $f$ factors through $\alpha$, but this provides us with a map $A \rightarrow X$, which factors through $f$. We obtain a map $\beta : T \rightarrow X$, and $\alpha\beta$ is an endomorphism of $T$ fulfilling $\alpha\beta f = f$, so invertible by minimality of the envelope. Thus $\alpha$ is a split epimorphism.
\ii Take $X \in \ST$ and any short exact sequence
\[
    \begin{tikzcd}
        0 \rar & X \rar & E \rar & T \rar & 0
    \end{tikzcd}
\]
$\ST$ is extension closed, so $E \in \ST$, and this extension splits since $T$ is a splitting projective.

Now consider \[\begin{tikzcd}
& A \dar[two heads] \drar["f"] &&& \\
0 \rar & \im f \rar["i"] & T \rar["c"] & \coker f \rar & 0\end{tikzcd}.\]
Applying $\Hom_A(-,X)$ we obtain
\[
\begin{tikzcd}[cramped, sep=scriptsize]
	\Hom_A(T,X) \rar["i^{*}"] \drar[two heads, "f^*", swap] & \Hom_A(\im f, X) \rar \dar["\cong"] & \Ext_A^{1}(\coker f, X) \rar & \stackrel{=0}{\Ext_A^{1}(T, X)} \\
	& \Hom_A(A, X) 
\end{tikzcd}
\]
The shown composition $\Hom_A(T,X) \rightarrow \Hom_A(A,X)$ has to be surjective since $f$ is a $\ST$-envelope. $\Hom_A(\im f, X) \rightarrow \Hom_A(A,X)$ is injective by left-exactness of $\Hom_A(-,X)$, and surjective since $f^*$ is, so an isomorphism. It follows that $i^{*}$ is surjective. Thus $\Ext_A^{1}(\coker f, X) = 0$.
\ii Take an $\add T'$ cover $\phi : W \twoheadrightarrow X$ of $X$, whose surjectivity is clear since $X \in \gen T'$. We obtain a short exact sequence
\[
\begin{tikzcd}
	0 \rar & \ker \phi \rar["\theta"] & W \rar["\phi"] & X \rar & 0
\end{tikzcd}
\]
Consider the map $r: A \rightarrow \ker \phi$, $a \mapsto au$ for some $u \in \ker \phi$. Since $T$ is a $\ST$-envelope the map $\theta r$ factors as $p f$ for some $p : T \rightarrow W$, and $\phi p$ factors as $q c$, where $q : \coker f \rightarrow X$.
\[
\begin{tikzcd}
	& A \rar["f"] \dar["r"] & T \rar["c"] \dar["p"] \dlar[dashed, "l", swap] & \coker f \rar \dar["q"] \dlar[dashed, "h", swap] & 0 \\
	0 \rar & \ker \phi \rar["\theta"] & W \rar["\phi"] & X \rar & 0
\end{tikzcd}
\]
The approximation property of $\phi$ allows us to find a map $h : \coker f \rightarrow W$, so that $\phi h = q$. Then $\phi (p - hc) = qc - qc = 0$, so we obtain a map $l : T \rightarrow \ker \phi$ with the property that $p - hc = \theta l$. 

But then $\theta r = pf = \theta l f + hcf = \theta l f$, since $cf = 0$, and since $\theta$ is a monomorphism we obtain that $r = l f$. Thus we obtain that our (arbitrarily chosen) element $u \in \ker \phi$ is contained in $\im(l)$. We repeat this process for a basis $u_1, \dots, u_m$ of $\ker \phi$, and obtain maps $l_1, \dots, l_m$ which assemble to a surjection $T^m \twoheadrightarrow \ker \phi$, proving the assertion.
\ii Let $S$ be an Ext-projective in $\ST$. We find an epimorphism $T'^n \twoheadrightarrow S$ whose kernel is in $\ST$ by (4). This sequence splits, so we're done.
\end{enumerate}
\end{proof}
\end{theorem}

\begin{remark}
    We will later show that this result—using some slightly more nuanced definitions, and tools from tilting theory—can be refined to a bijection between the functorially finite torsion classes $\ST$, and certain maximal \ta-rigid modules $T$, in the sense that adding any additional summand to $T$ would either generate a bigger torsion class, or lose \ta-rigidity. These are the so-called support \ta-tilting modules.

    This means that there might exist proper summands of $\SP(\ST)$ which generate all of $\ST$\footnote{Consider \cref{3KroneckerFFTorsionClasses}, or \cref{3LinearA3WithRelationExample}.}, and we need to account for that by restricting to exactly those \ta-rigid modules which are of the form $\SP(\ST)$.
\end{remark}

\begin{corollary}\label{3ffTorsionClassesEquivalence}
Let $(\ST, \SF)$ be a torsion theory in $A\Mod$. Then the following conditions are equivalent:
    \begin{enumerate}
        \ii $\ST$ is functorially finite.
        \ii $\ST = \gen (\SP(\ST))$
        \ii $\ST = \gen X$ for some $X$ in $A\Mod$.
    \end{enumerate}
\begin{proof}
    $1) \implies 2) \implies 3):$
    Assume that $\ST$ is functorially finite, and let $T'$ be as in the previous theorem. Then $\ST = \gen T'$, and $T'$ is Ext-projective in $\ST$, so $T' \in \add \SP(\ST)$. It follows that $\ST = \gen T' \subseteq \gen \SP(\ST) \subseteq \ST$.

    $3) \implies 1):$
    We know that $\gen X$ is a torsion class, and thus contravariantly finite. In \cref{2genMCovariantlyFinite} we showed that any module class of the form $\gen X$ is covariantly finite.
\end{proof}
\end{corollary}

\begin{example}[Non-functorially finite torsion classes]\label{3NonFfTorsionClassesExample}
    The Kronecker-quiver $\begin{tikzcd}
    Q: 1 \ar[r, shift right] \ar[r, shift left] & 2 \end{tikzcd}$ also gives us examples of torsion classes which are \emph{not} functorially finite. For this we define $\ST$ as the preinjective component of $KQ\Mod$, i.e. it contains all sums of modules of the form $\tau^{k}I$ for arbitrary $k$ in $\NN$, and arbitrary injective modules $I$. This is closed under extensions and quotients, making it a torsion class.

    By \cref{3ffTorsionClassesEquivalence} we know that, if $\ST$ were a functorially finite torsion class, then we would be able to find an Ext-projective module in $\ST$ which generates $\ST$.

    Take an arbitrary indecomposable $\tau^n I$ in $\ST$, where $n \in \NN$ and $I$ is an indecomposable injective module, and check whether it is Ext-projective. This cannot be the case since $\tau \tau^n I = \tau^{n+1} I$ is in $\ST$, which would contradict \cref{3extProjectivesBijection}. It follows that $\SP(\ST) = 0$, and thus $\ST$ is not functorially finite.

    Taking $\ST$ and adding all modules of an arbitrary set of regular "tubes" of the module category $KQ\Mod$ to $\ST$ provides further examples of non-functorially finite torsion classes over $KQ$. Quite a few of those, even, exactly one for each subset of $K \cup \{ \infty \}$, the set indexing the tubes.

    We recommend \cite{zbMATH07431617} for further insights into the lattice theory of torsion classes, including a summary of the torsion classes over the Kronecker-quiver.
\end{example}

\section{On the Number of Summands}

Recall that we write $\# X$ for the number of isomorphism classes of indecomposable summands of $X$, and that $\# A$ is the number of vertices of $A$.

We remind the reader of the definition of the support rank.

\begin{definition}[Support Rank]
    Given some $K$-algebra $A$, let $\SP(A\Mod)$ be a complete set of representatives of the isomorphism classes of its indecomposable projective modules. Given some module $X$ in $A\Mod$, we define the \emph{support rank} of $X$ as:
    \[
        \supportrank X := | \{ P \in \SP(A\Mod) \mid \Hom_A(P, X) \neq 0\}|,
    \]
    or in other words, as the number of isomorphism classes of indecomposable projective modules which have a non-zero map into $X$. From a quiver-theoretic point of view this is exactly the same as the number of vertices of $A$ at which $X$ is non-zero.

    Similarly, the support rank of a module class $\SC$ is defined as 
    \[
        \supportrank \SC := | \{ P \in \SP(A\Mod) \mid \exists C \in \SC : \Hom_A(P, C) \neq 0\}|.
    \]
    We call a module or a module class sincere if its support rank is maximal. 

    The support rank can just as easily be defined using injective modules, or by asking how many distinct simple modules arise in the composition series of a module, and these definitions are equivalent.
\end{definition}

With this the central result of this section is easily summarized: 
\begin{enumerate}
    \ii Every \ta-rigid module $T$ fulfills
        \[
            \# T \leq \# A,
        \]
    \ii If $T = \SP(\ST)$ for some functorially finite torsion class $\ST$, then 
        \[ 
            \# T = \supportrank T.
        \]
\end{enumerate}

We will see the importance of these results in the next section, when we focus our attention on precisely those modules for which this equality is given, the so-called support \ta-tilting modules.

\subsection{The Upper Bound}

\begin{definition}[Grothendieck Group]
    For any abelian category $\SA$, we define its Grothendieck group $K_0(\SA)$ to be the group generated by the (formal) symbols
    \[
        [X] 
    \]
    where $X$ ranges over all modules in $\SA$, modulo the relations
    \[
        [Y] = [X] + [Z]
    \]
    for any short exact sequence $0 \rightarrow X \rightarrow Y \rightarrow Z \rightarrow 0$ in $\SA$.

    Given an additive subcategory $\SC$ of some abelian category $\SA$, we define $K_0(\SC)$ in the same way, restricting ourselves to those short exact sequences fully contained in $\SC$. 

    For our purposes we merely need $K_0(A) := K_0(A\Mod)$, as well as $K_0(\proj A)$. Since all short exact sequences in $\proj A$ split, this is equivalent to considering the free $\ZZ$-module whose basis is given by the indecomposable projective modules.
\end{definition}

\begin{lemma}\label{4grothendieckGroupHasRankN}
    Given some $K$-algebra $A$, we see that $K_0(A)$ is the free $\ZZ$-module with basis given by the symbols of type $[S]$, where $S$ is a simple $A$-module.
    \begin{proof}
        Every module in $A\Mod$ has a finite composition series, which is essentially unique by the Jordan-Hölder Theorem.
    \end{proof}
\end{lemma}

To prove that every \ta-rigid module has at most $\# A$ summands, we first note that the \emph{dimension vector} $\ul{\dim}~M$ of any module $M$—using the above lemma—can be viewed as an element of $K_0(A)$.

\begin{definition}
    If $M$ is an $A$-module with minimal projective presentation
    \[
        \begin{tikzcd}
            P_1 \rar & P_0 \rar & M \rar & 0
        \end{tikzcd},
    \]
    then we define its \emph{g-vector} as $g(M) :=  \left[ P_0 \right] - \left[ P_1 \right]$, an element of $K_0(\proj A)$.

    We use this to define a bracket 
    \[
        \langle -, - \rangle : K_0(\proj A) \times K_0(A\Mod) \longrightarrow \ZZ ,
    \]
    given by 
    $
    \langle \left[P\right] , \ul{\dim}~M \rangle := \dim_K \Hom_A (P, M).
    $
\end{definition}

\begin{proposition}
	If $M, N$ are f.d. modules, then
    \[
        \langle g(M), \ul{\dim}~N \rangle = \dim_K \Hom_A(M, N) - \dim_K \Hom_A(N,\tau M)
    \]
\begin{proof}
    Namely, let $M$ have minimal projective presentation
    \[
        \begin{tikzcd}
            P_1 \rar & P_0 \rar & M \rar & 0
        \end{tikzcd},
    \]
    and apply $(-)^\vee$, then $D(-)$ to construct the sequence used to define the Auslander-Reiten translate $\tau M$. We apply $\Hom_A(N, -)$, and obtain 
    \[
        \begin{tikzcd}[column sep=1.05em]
            0 \rar & \Hom_A(N, \tau M) \rar &  \Hom_A(N, D(P_1^\vee)) \rar & \Hom_A(N,D(P_0^\vee))
        \end{tikzcd}
    \]
    We know that 
    \[
        \Hom_A(N,D(P^\vee)) \cong D(P^\vee \otimes_A N) \cong D \Hom_A(P, N),
    \]
    where the first isomorphism follows from Tensor-Hom-Adjunction. As for the second one: Any indecomposable projective module can be written as $Ae$, where $e$ is some idempotent, and the result then follows from $\Hom_A(Ae, A) \otimes_A N \cong eA \otimes_A N \cong eN \cong \Hom_A(Ae, N)$. We therefore obtain
    \[
        \begin{tikzcd}
            \Hom_A(P_0, N) \rar & \Hom_A(P_1, N) \rar & D \Hom_A(N, \tau M) \rar & 0
        \end{tikzcd}.
    \] 
    To complete the proof, we apply $\Hom_A(-, N)$ to the minimal projective presentation of $M$, which allows us to complete this exact sequence:
    \[
        \begin{tikzcd}[column sep=tiny]
            0 \rar & \Hom_A(M, N) \rar & \Hom_A(P_0, N) \rar & \Hom_A(P_1, N) \rar & D \Hom_A(N, \tau M) \rar & 0
        \end{tikzcd}
    \]
    The alternating sum of the dimensions is zero, which proves the claim.
\end{proof}
\end{proposition}

\begin{proposition}\label{4gVectorsDetermineTauRigid}
	Let $M, N$ be \ta-rigid modules fulfilling $g(M) = g(N)$. Then $M \cong N$.
\begin{proof}
	First prove that $M$ is generated by $N$. Since $N$ is \ta-rigid, $\gen N$ is a torsion class by \cref{3tauRigidModulesGenerateTorsionClasses}. Let $t_N M$ be the torsion submodule of $M$ relative to the torsion theory $(\gen N, N^{\perp 0})$. It is the sum of the images of all maps from $N$ to $M$. We get a canonical short exact sequence
    \[
    \begin{tikzcd}
        0 \rar & t_N M \rar & M \rar & F \rar & 0
    \end{tikzcd}
    \]
    where $F$ is torsion-free, so $\Hom_A(N, F) = 0$. We apply $\Hom_A(-, \tau M)$.
    \[
    \begin{tikzcd}
        0 \rar & \Hom_A(F, \tau M) \rar & \Hom_A(M, \tau M) \rar & \Hom_A(t_N M, \tau M)
    \end{tikzcd}
    \]
    Since $M$ is \ta-rigid, we get that $\Hom_A(F, \tau M) = 0$. 

    Now consider that 
    \[ 
    \begin{array}{cclc}
        \dim_K \Hom_A(M, F) &=& \dim_K \Hom_A(M, F) - \dim_K \Hom_A(F, \tau M) \\
                    &=& \langle g(M), \ul{\dim}~F \rangle \\ 
					&=& \langle g(N), \ul{\dim}~F \rangle \\
                    &=&  \dim_K \Hom_A(N, F) - \dim_K \Hom_A(F, \tau N) \\
                    &=& - \dim_K \Hom_A(F, \tau N).
    \end{array}
    \]
    This implies that $\dim_K \Hom_A(M, F) = 0$. But there is a map $M \twoheadrightarrow F$, so clearly $F = 0$. It follows that $t_N M = M$, so $M$ is in $\gen(N)$.

    By symmetry we obtain that $N \in \gen(M)$, and thus $\ST := \gen M = \gen N$. This is a torsion class, and both $M$ and $N$ are Ext-projective in $\ST$ by \cref{3tauRigidModulesGenerateTorsionClasses}. Write $\SF := \ST^{\perp 0}$.

    Note that this alone does not suffice to show that $\add M = \add N$, since a torsion class may be generated by several different sums of its Ext-projectives. We once again have to make use of $g(M) = g(N)$. For any $Y \in \SF$ we have $\Hom_A(M, Y) = 0 = \Hom_A(N, Y)$, and thus
    \[
        - \dim_K \Hom_A(Y, \tau M) = \langle g(M), \ul{\dim}~Y \rangle = \langle g(N), \ul{\dim}~Y \rangle = - \dim_K \Hom_A(Y, \tau N). \tag{$\ast$}
    \]
    Now take an $\add N$-cover $\phi : N_0 \rightarrow M$, which is surjective since $M \in \gen N$, and write $K := \ker \phi$. Applying $\Hom_A(N, -)$ to $0 \rightarrow K \rightarrow N_0 \rightarrow M \rightarrow 0$, and using that $\phi$ is an $\add N$-cover, we obtain a short exact sequence
    \[
        \begin{tikzcd}
            0 \rar & \Hom_A(N, K) \rar & \Hom_A(N, N_0) \rar & \Hom_A(N, M) \rar & 0.
        \end{tikzcd}
    \]
    Since $N_0, M \in \gen N$, we have $\Hom_A(N_0, \tau N) = 0 = \Hom_A(M, \tau N)$. Combining this with the additivity of $\langle g(N), - \rangle$ we get that
    \[
        \begin{array}{ccl}
            \langle g(N), \ul{\dim}~K \rangle &=& \langle g(N), \ul{\dim}~N_0 \rangle - \langle g(N), \ul{\dim}~M \rangle \\
                &=& \dim_K \Hom_A(N, N_0) - \dim_K \Hom_A(N, M) \\
                &=& \dim_K \Hom_A(N, K).
        \end{array}
    \]
    Comparing this with the previous proposition, which states that $\langle g(N), \ul{\dim}~K \rangle$ equals $\dim_K \Hom_A(N, K) - \dim_K \Hom_A(K, \tau N)$, we see that $\Hom_A(K, \tau N) = 0$. Consider the canonical short exact sequence $0 \rightarrow tK \rightarrow K \rightarrow K/tK \rightarrow 0$ of $K$ with respect to $(\ST, \SF)$. Since $K/tK$ is a quotient of $K$ we have $\Hom_A(K/tK, \tau N) = 0$, and since $K/tK \in \SF$, $(\ast)$ shows that $\Hom_A(K/tK, \tau M) = 0$. By the Auslander-Reiten formula, $\Ext_A^1(M, K/tK) = 0$. Furthermore $\Ext_A^1(M, tK) = 0$, since $tK \in \ST$ and $M$ is Ext-projective in $\ST$. From the exact sequence
    \[
        \begin{tikzcd}
            \Ext_A^1(M, tK) \rar & \Ext_A^1(M, K) \rar & \Ext_A^1(M, K/tK)
        \end{tikzcd}
    \]
    we obtain that $\Ext_A^1(M, K) = 0$. Thus $0 \rightarrow K \rightarrow N_0 \rightarrow M \rightarrow 0$ splits, and $M \in \add N$. By symmetry $N \in \add M$, so $\add M = \add N$.

    Thus $N$ and $M$ have the same indecomposable summands. Using the fact that $g(M) = g(N)$ and induction, we obtain that $M \cong N$: Take a common indecomposable summand $X$, and write $M = X \oplus M'$, $N = X \oplus N'$. Then $M'$ and $N'$ are \ta-rigid, and $g(M') = g(N')$, since $g$ is additive on direct sums.
\end{proof}
\end{proposition}

\begin{proposition}\label{4tauRigidUpperBound}
	Let $M$ be a \ta-rigid module with non-isomorphic indecomposable summands $M_1, \dots, M_n$. Then $g(M_1), \dots, g(M_n)$ are linearly independent over $\QQ$.

    In particular, any \ta-rigid module has at most $\# A$ distinct summands, i.e. $\# M \leq \# A$.
\begin{proof}
	Consider $\sum_{k=1}^{n} c_k g(M_k) = 0$ with $c_1, \dots, c_n \in \QQ$. Multiplying with an integer, we may assume that $c_1,\dots, c_n$ take on integer values. We're free to relabel our $M_k$ such that 
	\[
	\sum_{k=1}^{p} c_k g(M_k) = \sum_{k=p+1}^{p+q} (- c_k) g(M_k), \quad \text{ i.e. such that } \quad
	\begin{cases}
		c_{1}, \dots, c_{p} > 0 \\
		c_{p+1}, \dots, c_{p+q} < 0 \\
		c_{p+q+1}, \dots, c_{n} = 0
	\end{cases}
	\]
	Then 
	\[
	g(\bigoplus^p_{k=1} M_k^{c_k}) \quad = \quad g(\bigoplus^{p+q}_{k=p+1} M_k^{-c_{k}}).
	\]
	By the previous result, this shows that 
	\[
	\bigoplus^p_{k=1} M_k^{c_k} \quad \cong \quad \bigoplus^{p+q}_{k=p+1} M_k^{-c_{k}}
	\]
	Since the $M_k$ are pairwise non-isomorphic indecomposables, this is only possible if $p = q = 0$, i.e. if all $c_k$ are zero. This shows that $g(M_1), \dots,  g(M_n)$ are linearly independent. The upper bound on the number of summands of $M$ follows immediately from the rank of $K_0(\proj A)$.
\end{proof}
\end{proposition}

\subsection{Tilting Modules}
We now need some tools from tilting theory.

\begin{definition}
	We call a module $T$ a \emph{partial tilting module} if
	\begin{enumerate}
		\item $\Ext_A^1(T, T) = 0$
		\item $\projdim T \leq 1$
		
		We call a partial tilting module $T$ a \emph{tilting module} if the following is given.
		\ii $\codim_{T}(_{A}A) \leq 1$, i.e. there exists a short exact sequence
		\[
		\begin{tikzcd}
			0 \ar[r] & A \ar[r] & T^0 \ar[r] & T^1 \ar[r] & 0
		\end{tikzcd}
		\]
		with $T^0, T^1 \in \add T$.
	\end{enumerate}
    Note that every projective module is a partial tilting module, and that every projective generator is a tilting module. Also, note that every tilting module is faithful.
\end{definition}

\begin{theorem}[Bongartz' Completion]\label{4Bongartz}
	Let $M$ be a partial tilting module. There exists a module $E$ such that $T = E \oplus M$ is a tilting module. The module $E$ is constructed as follows: Take a basis $\xi_1,\cdots, \xi_n$ of $\Ext_A^1(M, A)$, and consider the tuple $(\xi_1, \cdots, \xi_n)$ as an element of $\Ext_A^1(M^n, A)$.
    The corresponding \emph{universal extension} defines $E$:
    \[
        \begin{tikzcd}
        0 \rar & A \rar & E \rar & M^n \rar & 0
        \end{tikzcd}
    \]
    \begin{proof}
        Apply $\Hom_A(M, -)$ to the universal extension and consider the long exact sequence: 
        \[
            \begin{tikzcd}[column sep=small]
                \Hom_A(M, M^n) \rar["\rho"] & \Ext_A^1(M, A) \rar & \Ext_A^1(M, E) \rar & \Ext_A^1(M, M^n)
            \end{tikzcd}
        \]
        Since $\Ext_A^1(M, M^n) = 0$ holds, $\Ext_A^1(M, E) = 0$ follows from the surjectivity of $\rho$. $\rho$ is surjective since $\rho$ maps the inclusion of $M$ into the $m$-th summand of $M^n$ onto $\xi_m$, and these are a basis. That $\Ext_A^1(E,M) = 0$ (resp. $\Ext_A^1(E, E) = 0$) follows with an application of $\Hom_A(-, M)$ (resp. $\Hom_A(-, E)$) to the universal extension, by noting that the surrounding modules in the long exact sequence are zero.
        Since $\projdim A, \projdim M \leq 1$, the same applies to $E$. Thus $M \oplus E$ is a partial tilting module. Finally, the universal extension $0 \rightarrow A \rightarrow E \rightarrow M^n \rightarrow 0$ is a short exact sequence of the form required in the definition of a tilting module, since $E, M^n \in \add(M \oplus E)$. This proves that $M \oplus E$ is a tilting module.
    \end{proof}
\end{theorem}

Once we have established further tools, we will be able to prove the existence of a Bongartz-style completion for \ta-rigid modules. We do this in \cref{5Bongartz}.

\begin{proposition}\label{4faithfulTauRigidIsPartialTilting}
\begin{enumerate}
    \ii[]
	\ii Any partial tilting module is \ta-rigid.
	\ii Any faithful \ta-rigid module is a partial tilting module.
\end{enumerate}
\begin{proof}
    Using \cref{1ProjDimEquiv}, the results follow easily.
    To prove the first claim, let $T$ be partial tilting. Since $\projdim T \leq 1$, no non-zero map from $T$ to $\tau T$ factors through an injective module, so $\Hom_A(T, \tau T) = \ol{\Hom}_A(T, \tau T) \cong D\Ext_A^1(T, T) = 0$.

    For the second claim, take a faithful \ta-rigid module $T$. That $\Ext_A^1(T, T) = 0$ is clear. We need to show that $\projdim T \leq 1$. Since $T$ is faithful, so is $DT$, and there exists an embedding $A\op \hookrightarrow DT^n$, thus $T^n \twoheadrightarrow DA$. Apply $\Hom_A(-, \tau T)$.
    \[
    \begin{tikzcd}
        0 \rar & \Hom_A(DA, \tau T) \rar & \Hom_A(T^n, \tau T) = 0
    \end{tikzcd}
    \]
    Thus $\Hom_A(DA, \tau T) = 0$, and so $\projdim T \leq 1$.
\end{proof}
\end{proposition}

Since any partial tilting module $T$ is \ta-rigid, it generates a torsion theory, which we will denote
\[
    (\ST, \SF) = (\gen T, T^{\perp 0}).
\]
We have the following chain of inclusions
\[
    \gen_1 T \subseteq \gen T \subseteq T^{\perp 1}.
\]
Surprisingly, tilting modules can be characterised precisely as those partial tilting modules for which these inclusions are equalities.

\begin{proposition}\label{4tiltingModuleEquivalentConditions}
    For a partial tilting module $T$ the following conditions are equivalent:
    \begin{enumerate}
        \ii $T$ is a tilting module.
        \ii $T^{\perp 0, 1} = 0$.
        \ii $\gen T = T^{\perp 1}$.
        \ii $\gen_1 T = T^{\perp 1}$.
        \ii $X$ is Ext-projective in $T^{\perp 1} \iff X \in \add T$.
    \end{enumerate}
    \begin{proof}
        All of these proofs boil down to applying a Hom-functor to the appropriate exact sequence, and noting that certain terms come out as zero.

        $1) \implies 2)$: Given some $X \in T^{\perp 0, 1}$, apply $\Hom_A(-, X)$ to the exact sequence necessary for $T$ to be a tilting module, and obtain 
        \[
            \begin{tikzcd}[row sep=0.0em]
            0 = \Hom_A(T^0, X) \rar & \Hom_A(A, X) \rar & \Ext_A^1(T^1, X) = 0
            \end{tikzcd}
        \]
        $2) \implies 3)$: $\gen T \subseteq T^{\perp 1}$ is clear. Take some $X$ so that $\Ext_A^1(T, X) = 0$. Our torsion theory provides a universal sequence
        \[
            \begin{tikzcd}
                0 \rar & tX \rar & X \rar & X/tX \rar & 0
            \end{tikzcd}
        \]
        and we apply $\Hom_A(T, -)$ to obtain 
        \[
            \begin{tikzcd}
                \Ext_A^1(T, X) \rar & \Ext_A^1(T, X/tX) \rar & \Ext_A^2(T, tX)
            \end{tikzcd}.
        \]
        The outer terms equal zero, the right one since $\projdim T \leq 1$, and thus so does the inner one. From this follows that $X/tX \in T^{\perp 0, 1}$, in other words, $X = tX \in \gen T$.

        $3) \implies 4)$: Take a $\add T$-(pre)cover $T' \xrightarrow{f} X$ of some $X \in T^{\perp 1}$, which is surjective since $X \in \gen T$. We consider the short exact sequence 
        \[
            \begin{tikzcd}
                0 \rar & \ker f \rar & T' \rar["f"] & X \rar & 0
            \end{tikzcd}
        \]
        and apply $\Hom_A(T,-)$ to obtain
        \[
            \begin{tikzcd}
                \Hom_A(T, T') \rar["f_*"] & \Hom_A(T,X) \rar & \Ext_A^1(T, \ker f) \rar & 0
            \end{tikzcd}
        \]
        in the associated long exact sequence. Since $f$ is an $\add T$-cover, $f_*$ is surjective, and $\Ext_A^1(T, \ker f) = 0$. Thus $\ker f$ is in $\gen T$, meaning that we can find an $\add T$-presentation for any $X$ in $T^{\perp 1}$.

        $4) \implies 5)$: Given that $X \in T^{\perp 1} = \gen_1 T$, we get to write a short exact sequence 
        \[
            \begin{tikzcd}
                0 \rar & \ker f \rar & T' \rar & X \rar & 0
            \end{tikzcd}
        \]
        where $\ker f \in \gen T$ and $T' \in \add T$, and thus $\Ext_A^1(X, \ker f) = 0$. It follows that this sequence splits, so $X \in \add T$. The other direction is clear.

        $5) \implies 1)$: Consider Bongartz' short exact sequence:
        \[
            \begin{tikzcd}
                0 \rar & A \rar & E \rar & T^n \rar & 0
            \end{tikzcd}
        \]
        We know that $E \oplus T$ is a tilting module. Let $X$ be a module such that $\Ext_A^1(T, X) = 0$. Apply $\Hom_A(-, X)$ to the above sequence:
        \[
            \begin{tikzcd}
                \Ext_A^1(T^n, X) \rar & \Ext_A^1(E, X) \rar & \Ext_A^1(A, X)
            \end{tikzcd}
        \]
        The outer terms are zero, and so is the inner one. Since $E \in T^{\perp 1}$, it follows that $E$ is Ext-projective in $T^{\perp 1}$. It follows that $E \in \add T$, so the above sequence shows that $T$ is already a tilting module.
    \end{proof}

\end{proposition}

The following result makes handling of tilting modules substantially easier, and is a key ingredient to develop the theory of \ta-tilting modules.

\begin{proposition}\label{4partialTiltingBound}
    Any partial tilting module ${}_A T$ fulfills 
    \[
        \#T \leq \#A,
    \]
    with equality if and only if ${}_A T$ is a tilting module.
    \begin{proof}
        Let ${}_A T$ be a partial tilting module. By \cref{4faithfulTauRigidIsPartialTilting} we know that $T$ is a \ta-rigid module. \Cref{4tauRigidUpperBound} then shows $\# T \leq \# A$.

        Now let $T$ be a tilting module. We first show that there exists a short exact sequence
       \[
            \begin{tikzcd}
                0 \ar[r] & P \ar[r, "f"] & T^0 \ar[r] & T^1 \ar[r] & 0
            \end{tikzcd}
		\]
        for every indecomposable projective module $P$, where $T^0$ and $T^1$ are in $\add T$. For this, take an $\add T$-preenvelope $f : P \rightarrow T^0$. By the definition of a tilting module, $P \leq A$ embeds into a module in $\add T$, and this embedding factors through $f$, so $f$ is injective. Set $T^1 := \coker f$, which is contained in $\gen T = T^{\perp 1}$ by \cref{4tiltingModuleEquivalentConditions}. Given any $X \in T^{\perp 1} = \gen T$, we apply $\Hom_A(-, X)$ and obtain the exact sequence
        \[
            \begin{tikzcd}
                \Hom_A(T^0, X) \rar["f^*"] & \Hom_A(P, X) \rar & \Ext_A^1(T^1, X) \rar & \Ext_A^1(T^0, X) = 0.
            \end{tikzcd}
        \]
        Here $f^*$ is surjective: Any map $P \rightarrow X$ lifts along a surjection $T^n \twoheadrightarrow X$ since $P$ is projective, and this lift factors through $f$. Thus $\Ext_A^1(T^1, X) = 0$, so $T^1$ is Ext-projective in $T^{\perp 1}$, and $T^1 \in \add T$ by \cref{4tiltingModuleEquivalentConditions}.

        We now work in $K_0(\proj A)$. All modules in $\add T$ have projective dimension at most one, so their minimal projective presentations are projective resolutions $0 \rightarrow Q_1 \rightarrow Q_0 \rightarrow X \rightarrow 0$, and $g(X) = [Q_0] - [Q_1]$. This formula holds for \emph{any} projective resolution of length one, since such a resolution is the direct sum of the minimal one and a complex of the form $Q \xrightarrow{\id} Q$. Applying the Horseshoe lemma to the above short exact sequence we obtain a projective resolution of $T^0$ of length one, and thus
        \[
            \left[ P \right] = g(T^0) - g(T^1).
        \]
        Since the indecomposable projective modules form a basis of $K_0(\proj A) \cong \ZZ^{\# A}$, the $g$-vectors of the indecomposable summands of $T$ generate a group of rank $\# A$. Consequently, $\# T \geq \# A$, and every tilting module has exactly $\# A$ distinct summands.

        We warn the reader that it would not suffice to argue with the classes $[T^0] - [T^1]$ in $K_0(A)$: The indecomposable projective modules need not generate a group of rank $\# A$ in $K_0(A)$ if $A$ has infinite global dimension. For example, if $A$ is given by the quiver with two vertices and arrows $1 \rightleftarrows 2$ modulo all paths of length two, then $[P_1] = [S_1] + [S_2] = [P_2]$.

        What remains to be shown is that any partial tilting module with $\# A$ summands is already a tilting module. To see this, we take a partial tilting module $T$ which fulfills $\# T = \# A$, and use Bongartz' completion (\cref{4Bongartz}) to complete $T$ to a tilting module $T \oplus E$. 

        Of course, $\# A = \# T \leq \#(T \oplus E) \leq \# A$. It follows that $E$ is an element of $\add T$, and thus the universal extension $0 \rightarrow A \rightarrow E \rightarrow T^n \rightarrow 0$ shows that $T$ is already a tilting module. This completes the proof.
    \end{proof}
\end{proposition}

\section{\texorpdfstring{$\tau$}{Tau}-tilting Modules}
\subsection{Support Rank}
In the following we return to our original goal, the study of \ta-rigid modules and their associated torsion classes. We use what we learned about tilting modules—first and foremost \cref{4partialTiltingBound}—to establish similar results for \ta-rigid modules, and to develop a more nuanced understanding of their overall structure.

First, an important lemma.
\begin{lemma}\label{4tauPreservedWhenPassingToQuotientAlgebras}
    Let $A$ be a finite dimensional algebra, and $I$ an ideal in $A$. Let $M, N$ be $(A/I)$-modules.

    \begin{enumerate}
        \ii $\Hom_A(N, \tau M) = 0 \quad \implies \quad \Hom_{A/I}(N, \tau_{A/I} M) = 0$.
    \end{enumerate}
    If $I = \langle e \rangle$ for an idempotent $e$ in $A$ holds, then we even have:
    \begin{enumerate}
        \ii[2)] $\Hom_A(N, \tau M) = 0 \quad \iff \quad \Hom_{A/I}(N, \tau_{A/I} M) = 0$.
    \end{enumerate}
    \begin{proof}
        We have a natural inclusion $\Ext_{A/I}^1(M,N) \hookrightarrow \Ext_A^1(M,N)$. If $I = \langle e \rangle$, then the additivity of $\dim_K$ at each vertex ensures that this is an isomorphism.
        Assume that $\Hom_A(N, \tau M) = 0$ holds. Using the Auslander-Smal\o\ lemma, the above observation, and the lemma again, we obtain that:
        \begin{align*}
\Hom_A(N, \tau M) = 0 &\iff \Ext^1_A(M, \gen N) = 0 \\
                      & \implies  \Ext^1_{A/I}(M, \gen N) = 0 & \\
                      & \iff \Hom_{A/I}(N, \tau_{A/I} M) = 0
        \end{align*}
        The proof of the second result is analogous.
    \end{proof}
\end{lemma}

\begin{proposition}\label{4PTIsTiltingOverAAnn}
Let $(\ST, \SF)$ be a torsion theory in $A\Mod$. Then the following conditions are equivalent:
    \begin{enumerate}
        \ii $\ST$ is functorially finite.
        \ii $\ST = \gen (\SP(\ST))$
        \ii $\ST = \gen X$ for some $X$ in $A\Mod$.
        \ii $\SP(\ST)$ is a tilting $(A/\ann \ST)$-module.
    \end{enumerate}
\begin{proof}
    The equivalence of the first three statements was already proved in \cref{3ffTorsionClassesEquivalence}. Prove that those three imply $(4)$. For this, write $T := \SP(\ST)$, as well as $A_0 := (A/\ann \ST)$. First of all, note that $A_0\Mod$ embeds into $A\Mod$, and that all modules in $\ST$ can interchangeably be viewed as $A_0$-modules, or as $A$-modules.

    $1),~2),~3) \implies 4):$ Prove that $T$ is faithful as an $A_0$-module. Take an element $a \in A_0$ such that $aT = 0$. By \cref{AuslanderSmaloTheorem} we know that every module $X$ in $\ST$ can be written as a quotient of $T^n$, for some $n$. It follows that $aX = 0$ for all modules $X$ in $\ST$, so $a$ has to equal zero.

    By the previous \cref{4tauPreservedWhenPassingToQuotientAlgebras} we know that $T$ is \ta-rigid over $A_0$. Thus by \cref{4faithfulTauRigidIsPartialTilting} we obtain that $T$ is in fact a partial tilting module over $A_0$. It remains to be shown that $T$ is a tilting module. 

    Let $f_A : A \rightarrow T^0$ be a $\ST$-envelope of $A$, and write $T^1 := \coker f_A$. By \cref{AuslanderSmaloTheorem}, $T^0$ and $T^1$ are Ext-projective in $\ST$, and thus contained in $\add T$. Since $T^0 \in \ST$ is annihilated by $\ann \ST$, the map $f_A$ factors as $A \twoheadrightarrow A_0 \xrightarrow{f} T^0$. The map $f$ is a $\ST$-envelope of $A_0$: Any map from $A_0$ to a module in $\ST$ factors through $f_A$, and thus through $f$, since $A \twoheadrightarrow A_0$ is an epimorphism. Similarly, if $\alpha \in \End_A(T^0)$ fulfills $\alpha f = f$, then $\alpha f_A = f_A$, so $\alpha$ is invertible. Since $A \twoheadrightarrow A_0$ is surjective, we also have $\coker f = \coker f_A = T^1$, and we obtain an exact sequence
    \[
        \begin{tikzcd}
            A_0 \rar["f"] & T^0 \rar & T^1 \rar & 0
        \end{tikzcd}
    \]
    with $T^0, T^1 \in \add T$. It remains to show that $f$ is injective. We have shown above that $T$ is a faithful $A_0$-module. Choosing a $K$-basis $m_1, \dots, m_r$ of $T$, the map $A_0 \rightarrow T^r$, $a \mapsto (a m_1, \dots, a m_r)$ is injective, since its kernel annihilates $T$. As $T^r \in \ST$, this map factors through the $\ST$-envelope $f$, so $f$ is injective as well. This proves that $T$ is a tilting module.

    $4) \implies 3):$ We know that $T$ is a tilting module over $A_0$. By \cref{4tiltingModuleEquivalentConditions} we see that $\gen T = T^{\perp 1}$, where $T^{\perp 1}$ is taken in $A_0\Mod$. Given $X \in \ST$, we have $\Ext^1_{A_0}(T, X) \subseteq \Ext^1_A(T, X) = 0$, since $T$ is Ext-projective in $\ST$. It follows that $\ST \subseteq T^{\perp 1} = \gen T$. Of course, $\gen T \subseteq \ST$. The claim follows.
\end{proof}
\end{proposition}

\begin{proposition}\label{4torsionClassSupportRankBound}
    For any functorially finite torsion class $\ST$ we obtain that
    \[
        \#\SP(\ST) = \supportrank(\ST).
    \]
    \begin{proof}
        Take any functorially finite torsion class $\ST$. By \cref{4PTIsTiltingOverAAnn}, $\SP(\ST)$ is a tilting $(A/\ann\ST)$-module. In particular, 
        \[
            \#\SP(\ST) = \#(A/\ann\ST) = \supportrank(\SP(\ST)) = \supportrank(\ST).
        \]
        The first equality is an immediate consequence of \cref{4partialTiltingBound}, the second follows since $A/\ann \ST$ embeds into $\SP(\ST)$, and the third since $\SP(\ST)$ generates $\ST$.
    \end{proof}
\end{proposition}

\subsection{Support \texorpdfstring{$\tau$}{Tau}-tilting Modules}
\begin{definition}[\ta-tilting modules]
    Let $T$ be a \ta-rigid $A$-module.
    \begin{enumerate}
        \ii We call $T$ a \emph{\ta-tilting module} if $\# A = \# T$.
        \ii We call $T$ a \emph{support \ta-tilting module} if there exists an idempotent $e$ of $A$ such that $T$ is a \ta-tilting $(A/\langle e \rangle)$-module.
    \end{enumerate}
\end{definition}

\begin{proposition}\label{5tauRigidFullRankIsSupportTauTilting}
    A \ta-rigid module $T$ fulfills
    \[
        \# T \leq \supportrank T,
    \]
    with equality if and only if $T$ is a support \ta-tilting module. As a consequence, \ta-tilting modules can be precisely characterised as sincere support \ta-tilting modules.
    \begin{proof}
        Since $T$ is \ta-rigid, we consider the associated torsion class $\ST = \gen T$. Now consider $\SP(\ST)$. \Cref{4torsionClassSupportRankBound} gives us that
        \[
            \# T \leq \# \SP(\ST) = \supportrank(\ST) = \supportrank(T).
        \]
        Now assume $T$ is support \ta-tilting. In that case
        \[
            \supportrank T \leq \# (A/\langle e \rangle) = \#T \leq \supportrank T.
        \]
        These (in)equalities follow since $T$ is an $(A/\langle e \rangle)$-module, since $T$ is \ta-tilting over $(A/\langle e \rangle)$, and by the above.
        
        We now merely need to prove the other direction. Assume that $T$ is a \ta-rigid module fulfilling $\# T = \supportrank T$. Take a maximal basic projective module $P := Ae$, where $e$ is an idempotent of $A$ fulfilling $\Hom_A(P, T) = 0$. By \cref{4tauPreservedWhenPassingToQuotientAlgebras} $T$ is a \ta-rigid $(A/\langle e \rangle)$-module, and we get that 
        \[
            \# (A/\langle e \rangle) = \supportrank T = \# T,
        \]
        which shows that $T$ is support \ta-tilting.
    \end{proof}
\end{proposition}

We shine light on the relationship between \ta-tilting modules and tilting modules:

\begin{proposition}\label{5tiltingAreFaithfulTauTilting}
    \begin{enumerate}
        \ii[]
        \ii Tilting modules are precisely faithful \ta-tilting modules.
        \ii Any \ta-rigid module $T$ is a partial tilting $(A/\ann T)$-module.
        \ii Any \ta-tilting module $T$ is a tilting $(A/\ann T)$-module.
    \end{enumerate}
    \begin{proof}
        The reader is encouraged to look at \cref{4faithfulTauRigidIsPartialTilting}, from which we obtain that a faithful \ta-rigid module is partial tilting. It is clear that a \ta-rigid module $T$ is faithful over $(A/\ann T)$, and by \cref{4tauPreservedWhenPassingToQuotientAlgebras} it is also \ta-rigid over $(A/\ann T)$. This already shows $2)$.

        Using \cref{4partialTiltingBound} and the fact that $T$ has $\# A$ summands, we obtain $3)$, and one direction of $1)$.

        To complete the proof of $1)$, note that—also by \cref{4faithfulTauRigidIsPartialTilting}—any tilting module is \ta-rigid. Using the short exact sequence in the definition of a tilting module we see that any tilting module is necessarily faithful.
    \end{proof}
\end{proposition}

A convenient view of support \ta-tilting modules is as pairs of modules with certain properties. This view is fully compatible, as we will show momentarily, and it allows us to make important information related to the support rank of \ta-rigid modules more explicit.

\begin{definition}[\ta-tilting pairs]
    Let $T$ be an $A$-module, and let $P$ be a projective module.
    \begin{enumerate}
        \ii We call $(T,P)$ a \emph{\ta-rigid pair} if $T$ is \ta-rigid, and $\Hom_A(P, T) = 0$.
        \ii We call $(T,P)$ a \emph{support \ta-tilting pair} if $(T,P)$ is a \ta-rigid pair, and $\#T + \#P = \# A$.
    \end{enumerate}
    We call such a pair \emph{basic}, if $T$ and $P$ are basic, or \emph{multiplicity-free} modules, i.e. any two distinct indecomposable summands of $T$ (resp. of $P$) are not isomorphic to each other.
\end{definition}

\begin{proposition}\label{5tauRigidPairsAreModules}
    Let $(T, P)$ be a pair of $A$-modules, such that $P$ is a projective module. Let $e$ be an idempotent of $A$ such that $\add P = \add Ae$. 
    \begin{enumerate}
        \ii $(T,P)$ is a \ta-rigid pair for $A$ if and only if $T$ is a \ta-rigid $A/\langle e \rangle$ module.
        \ii $(T,P)$ is a support \ta-tilting pair for $A$ if and only if $T$ is a \ta-tilting $(A/\langle e \rangle)$-module.
        \ii If $(T,P)$ and $(T, P')$ are support \ta-tilting pairs for $A$, then $\add P = \add P'$. In other words, $T$ determines $P$ up to isomorphism.
    \end{enumerate}
    \begin{proof}
        \begin{enumerate}
            \ii Clear, using \cref{4tauPreservedWhenPassingToQuotientAlgebras}.
            \ii Note that $\#(A/\langle e \rangle) = \# A - \# P = \# T$.
            \ii $T$ is a support \ta-tilting module, and thus fulfills $\# T = \supportrank T$. Using that $\#T + \#P = \#A$, we obtain that $\#P = \#A - \supportrank T$. There are exactly $\#A - \supportrank T$ isomorphism classes of indecomposable projective modules in ${}^{\perp 0} T$, and every indecomposable summand of $P$ is one of them. Thus $P$, and likewise $P'$, has every one of them as a summand, so $\add P = \add P'$.
        \end{enumerate}
    \end{proof}
\end{proposition}

We provide a few characterisations of what it means to be a support \ta-tilting module.
\begin{proposition} \label{5supportTauRigidEquivalences}
    For a basic \ta-rigid module $T$ the following conditions are equivalent:
    \begin{enumerate}
        \ii $T$ is a support \ta-tilting module.
        \ii $(T, P)$ is a support \ta-tilting pair, where $P$ is a maximal basic projective module fulfilling $\Hom_A(P, T) = 0$.
        \ii $T$ is a \ta-tilting $(A/\langle e \rangle)$-module, where $e$ is a maximal idempotent of $A$ such that $eT = 0$.
        \ii $\# T = \supportrank T$.
        \ii $\SP(\gen T) \cong T$.
        \ii $\SP(\ST) \cong T$  for some functorially finite torsion class $\ST$.
    \end{enumerate}
    \begin{proof}
        $1) \iff 4):$ This is \cref{5tauRigidFullRankIsSupportTauTilting}.

        $4) \implies 2):$ We take such a $P$. Naturally, $(T, P)$ is a \ta-rigid pair. We know that $\# T = \supportrank T = \# A - \# P$, so the conclusion follows.

        $2) \implies 3):$ Note that $P = Ae$. Using \cref{5tauRigidPairsAreModules}, $T$ is a \ta-rigid $(A/\langle e \rangle)$-module, and $T$ clearly fulfills $\# (A/\langle e \rangle) = \# A - \# P = \# T$, making it \ta-tilting.

        $3) \implies 5):$ We know that
        \[
            \# \SP(\gen T) = \supportrank (\gen T) = \# (A/\langle e \rangle) = \# T.
        \]
        These equalities follow using \cref{4torsionClassSupportRankBound}, using that $T$ is a sincere $(A/\langle e \rangle)$-module by maximality of $e$, and from the fact that $T$ is a \ta-tilting $(A/\langle e \rangle)$-module respectively. Since $T$ is Ext-projective in $\gen T$, we have $\add T \subseteq \add \SP(\gen T)$, so we obtain that $\add \SP(\gen T) = \add T$. Since $\SP(\gen T)$ and $T$ are basic, the result follows.

        $5) \implies 6):$ This is clear, since $\gen T$ is functorially finite by \cref{3tauRigidsGenerateFFTorsionClasses}.

        $6) \implies 4):$ Since $\ST$ is functorially finite, $\ST = \gen \SP(\ST) = \gen T$ by \cref{3ffTorsionClassesEquivalence}. Now the result follows from \cref{4torsionClassSupportRankBound}, since $\# T = \# \SP(\ST) = \supportrank \ST = \supportrank T$.

        The assumption that $\ST$ is functorially finite in $6)$ cannot be dropped. For example, let $A$ be the path algebra of $1 \rightarrow 2 \rightrightarrows 3$, and let $\ST$ be the smallest torsion class containing $P_1$ and all regular and preinjective modules over the Kronecker quiver $2 \rightrightarrows 3$. Then $\SP(\ST) = P_1$, but $\# P_1 = 1 < 3 = \supportrank P_1$.
    \end{proof}
\end{proposition}

\subsection{Bijections}
In this section we will prove some of the fundamental results that justify the concept of \ta-tilting modules.

For this, we use \ftors $A$ to refer to the set of functorially finite torsion classes in $A\Mod$, and we refer to the set of isomorphism classes of basic support \ta-tilting modules as \sttilt $A$.

First, we will make the intuition explicit that there is in fact a bijection between \ftors $A$, and \sttilt $A$.

Second, we will look at the relationship between \sttilt $A$ and \sttilt $A\op$, as well as at the relationship between \sttilt $A$, and an analogously defined class of support $\tau^-$ tilting modules. It turns out that we will be able to find natural bijections between these sets.

\begin{theorem}\label{5sTauTiltFTorsBijection}
There is a bijection \[ 
\begin{array}{c}
	\sttilt A \\
	T \\
	\SP(\ST)
\end{array} 
\begin{array}{c}
	\xleftrightarrow{\enspace \, \sim \,\enspace} \\
	\xrightarrow{\enspace \, \quad \,\enspace} \\
	\xleftarrow{\enspace \quad \enspace}
\end{array}
\begin{array}{c}
	\ftors A \\
	\gen T \\
	\ST 
\end{array}
\]

\begin{proof}
    Take a support \ta-tilting module $T$. Using \cref{3tauRigidsGenerateFFTorsionClasses}, we know that $\gen T$ is a functorially finite torsion class. By \cref{5supportTauRigidEquivalences} we know that $\SP(\gen T) \cong T$.

    Now take a functorially finite torsion class $\ST$. By \cref{3ffTorsionClassesEquivalence} we know that $\SP(\ST)$ is a basic \ta-rigid module which generates $\ST$. Using \cref{5supportTauRigidEquivalences} we know that $\SP(\ST)$ is a support \ta-tilting module.

    Consequently, these bijections are well-defined and mutually inverse.
\end{proof}
\end{theorem}

\begin{corollary}
    The bijection of the preceding theorem induces bijections
    \[
        \begin{tikzcd}
            \ttilt A \rar[leftrightarrow] & \sftors A
        \end{tikzcd}
        \qquad \text{and} \qquad
        \begin{tikzcd}
            \tilttilt A \rar[leftrightarrow] & \fftors A
        \end{tikzcd},
    \]
    where 
    \begin{itemize}
        \ii $\ttilt A$ is the class of isomorphism classes of basic \ta-tilting modules over $A$.
        \ii $\sftors A$ is the class of sincere functorially finite torsion classes over $A$.
        \ii $\tilttilt A$ is the class of isomorphism classes of basic tilting modules over $A$.
        \ii $\fftors A$ is the class of faithful functorially finite torsion classes over $A$.\footnote{Analogously to the definition of a faithful module, a module class $\SC$ is faithful if its annihilator, i.e. the ideal of all elements $a \in A$ which act as zero on all modules in $\SC$, is the zero ideal.}
    \end{itemize}
    \begin{proof}
        Using \cref{5tauRigidFullRankIsSupportTauTilting}, each \ta-tilting module is sincere. \Cref{4torsionClassSupportRankBound} gives us that $\SP(\ST)$ is a \ta-tilting module when $\ST$ is sincere.

        Similarly, since tilting modules are faithful—by \cref{5tiltingAreFaithfulTauTilting}—it is clear that they generate faithful torsion classes, and those necessarily have a faithful \ta-rigid generator.
    \end{proof}
\end{corollary}

By switching our focus from $\tau$ to $\tau^-$, we are free to think of $\tau^-$-rigid modules.

\begin{definition}
    We call an $A$-module $X$ $\tau^-$-rigid if
    \[
        \Hom_A(\tau^- X, X) = 0.
    \]
    We write $\cogen X$ for the modules cogenerated by $X$, i.e. those modules that inject into some module in $\add X$.
\end{definition}

\Cref{3ExampleSkewedTriangle} provides an example of a module which is \ta-rigid, but not $\tau^-$-rigid, and vice versa.

\begin{proposition}
    An $A$-module $X$ is \ta-rigid exactly when $DX$ is a $\tau^-$-rigid $A\op$-module.
    \begin{proof}
        As the Auslander-Reiten translations are defined as 
        \[
            \tau = D \circ \text{Tr} \quad \text{ and } \quad \tau^- = \text{Tr} \circ D
        \]
        we use that $\Hom_{A\op}(Y,Z) \cong \Hom_{A}(DZ,DY)$ to obtain
        \[
            \Hom_{A\op}(\tau^- DX, DX) \cong \Hom_{A}(DDX, D\text{Tr}DDX) \cong \Hom_{A}(X, \tau X).
        \]
    \end{proof}
\end{proposition}

\begin{remark}\label{5dualAuslanderSmaloLemma}
    The dual of the Auslander-Smal\o{} Lemma is the equivalence of the following two statements:
    \begin{enumerate}
        \ii $\Hom_A(\tau^- N, M) = 0$
        \ii $\Ext_A^1(\cogen M, N) = 0$
    \end{enumerate}
    Similarly, other results dualize, in particular \cref{2genMCovariantlyFinite}, \cref{3tauRigidModulesGenerateTorsionClasses}, and \cref{AuslanderSmaloTheorem}, as well as \cref{4PTIsTiltingOverAAnn}.

    For our purposes an explicit theory of support $\tau^-$-tilting modules would be redundant. We will now see why.
\end{remark}

\begin{definition}
    We write $\trigid A$ for the set of isomorphism classes of basic \ta-rigid pairs of $A$.

    Similarly, $\tminusrigid A$ is the set of isomorphism classes of pairs $(U,I)$, where $U$ is a basic $\tau^-$-rigid module, and $I$ is a basic injective module such that $\Hom_A(U, I) = 0$. Dualising \cref{3defExtProjectives}, we write $\SI(\SF)$ for a direct sum consisting of indecomposable Ext-injective modules in $\SF$, choosing one module for each isomorphism class, and where $\SF$ is a torsion-free class. 

    A support $\tau^-$-tilting pair is a $\tau^-$-rigid pair $(U,I)$ fulfilling $\# U + \# I = \# A$. The set of isomorphism classes of support $\tau^-$-tilting pairs will be denoted $\sminusttilt A$.
\end{definition}

\begin{proposition}\label{5definitionDagger}
    Take a \ta-rigid $A$-module $T$. We decompose $T$ as $T = T\dualpr \oplus T\dualnp$, where $T\dualpr$ is a maximal projective summand of $T$. For any \ta-rigid pair $(T, P)$ we define
    \[
        (T, P)^\dagger := (\transpose T\dualnp \oplus P^*, T\dualpr^*) = (\transpose T \oplus P^*, T\dualpr^*),
    \]
    where $(-)^* := \Hom_A(-, {}_AA)$. We obtain the following result:

    $(-)^\dagger$ gives bijections
    \[
        \trigid A \longleftrightarrow \trigid A\op \qquad \text{ and } \qquad \sttilt A \longleftrightarrow \sttilt A\op
    \]
    such that $(-)^{\dagger\dagger} = \id$.
    \begin{proof}
        Remembering that the transpose acts as zero on projective modules, and gives an involution on non-projective modules, is enough to show that $(-)^\dagger$ is an involution. We now merely have to show that $(T, P)^\dagger$ is a \ta-rigid pair for $A\op$. We start by showing that $\transpose T \oplus P^*$ is \ta-rigid. We use that $T$ is \ta-rigid, as well as duality for finite-dimensional vector spaces, to obtain
        \[
            0 = \Hom_A(T\dualnp, \tau T) = \Hom_{A\op}(\transpose T, D T\dualnp) = \Hom_{A\op}(\transpose T, \tau \transpose T).
        \]
        Noting that $P = Ae$ for some idempotent $e$, we see that $\Hom_A(Ae, A) \cong eA$, and $P^* \otimes_A T\dualnp = 0$ follows since $\Hom_A(P, T) = 0$. Using tensor-hom adjunction, we see that
        \[
            0 = D(P^* \otimes_A T\dualnp) = \Hom_{A\op}(P^*, DT\dualnp) = \Hom_{A\op}(P^*, \tau \transpose T).
        \]
        Combining the above, $\transpose T \oplus P^*$ is \ta-rigid. That there are no maps from $T\dualpr^*$ to $\transpose T$ follows from,
        \[
            0 = \Hom_A(T\dualpr, \tau T) = \Hom_{A\op}(\transpose T, DT\dualpr) = D\Hom_{A\op}(T\dualpr^*, \transpose T)
        \]
        using the same tools as above, as well as from $\Hom_{A\op}(T\dualpr^*, P^*) = 0$, which is given since $\Hom_A(P, T\dualpr) = 0$. This shows that $(T, P)^\dagger$ is a \ta-rigid pair.

        Finally, $\transpose T\dualnp$ has no projective summands, and $\transpose$ as well as $(-)^*$ preserve the number of non-isomorphic indecomposable summands. Thus $\#(\transpose T\dualnp \oplus P^*) + \# T\dualpr^* = \# T\dualnp + \# P + \# T\dualpr = \# T + \# P$, so $(-)^\dagger$ sends support \ta-tilting pairs to support \ta-tilting pairs. This completes the proof.
    \end{proof}
\end{proposition}

These bijections are interesting on their own terms. That there is a left-right symmetry for \ta-rigid modules in particular is surprising. Generally, similar results relating classes of support (co)tilting modules over some algebra $A$ and $A\op$ are not nearly as nice, if they exist at all. This is just one area in which \ta-tilting theory is substantially better behaved than classical tilting theory.

\begin{corollary}\label{5tauTauMinusBijection}
    Composing this map with standard vector space duality we obtain bijections
    \[
        (X, P) \longmapsto (\tau X \oplus \nu P, \nu X\dualpr)
    \]
    between $\trigid A$ and $\tminusrigid A$, and between $\sttilt A$ and $\sminusttilt A$. Here $X\dualpr$ is a maximal projective summand of $X$, as in \cref{5definitionDagger}, and dually $Y_{\text{in}}$ is a maximal injective summand of $Y$.
    \[
    \begin{array}{ccc}
        \trigid A & \xleftrightarrow{\quad \sim \quad} & \tminusrigid A \\
        (X, P) & \xrightarrow{D ((-)^\dagger)} & (\tau X \oplus \nu P, \nu X\dualpr) \\
        (\tau^- Y \oplus \nu^- I, \nu^- Y_{\text{in}}) & \xleftarrow{(D(-))^\dagger} & (Y, I)
    \end{array}
    \]
\end{corollary}

This allows us to answer a remarkably natural question the reader might have had. The proof is a little fiddly, but not very complicated.

\begin{theorem}\label{5STFfIffSFFf}
    Take a torsion theory $(\ST, \SF)$.
    
    $\ST$ is functorially finite if and only if $\SF$ is functorially finite.
    \begin{proof}
        We prove that $\ST$ is functorially finite if $\SF$ is. The converse is dual: Apply the statement to the torsion theory $(D\SF, D\ST)$ in $A\op\Mod$.

        Assume that $\SF$ is functorially finite. Recall from \cref{3defExtProjectives} that $\SP(\ST)$ is a module with at most $\# A$ summands. We write $T := \SP(\ST)$, and take a maximal basic projective module $P$ which fulfills that $\Hom_A(P, \ST)~=~0$. To prove that $\ST$ is functorially finite, we aim to use \cref{4PTIsTiltingOverAAnn}, and thus need to show that $T$ is a tilting module over $A' := (A/\ann \ST)$.

        First, note that we have embeddings $\ST \hookrightarrow A'\Mod \hookrightarrow A\Mod$.

        $\ST$ contains a module whose annihilator coincides with that of the torsion class. To see this, take a module $M \in \ST$ such that $\ann M$ has minimal $K$-dimension. Clearly $\ann \ST = \bigcap_{N \in \ST} \ann N \subseteq \ann M$. If $\dim \ann \ST < \dim \ann M$, then there exists some element $a \in \ann M$, and a module $N \in \ST$, such that $a N \neq 0$. Clearly, $\ann (M \oplus N) \subsetneq \ann M$. This would contradict the assumption that $\ann M$ has minimal $K$-dimension, so $\dim \ann \ST = \dim \ann M$, and the statement follows.

        By the above, we in particular obtain that there exists a faithful $A'$-module in $\ST$. An immediate consequence of this is that $\ST$ contains all injective $A'$-modules\footnote{Consider the dual: If $U$ is faithful, then $A$ injects into $U^n$ for some $n$, so $A \in \cogen U$.}. From this follows that $\Hom_{A}(DA', \tau T) = 0$, since $\tau T \in \SF$ by \cref{3extProjectivesBijection}. By \cref{4tauPreservedWhenPassingToQuotientAlgebras} we get $\Hom_{A'}(DA', \tau_{A'} T) = 0$. Using \cref{1ProjDimEquiv}, we obtain that $\projdim_{A'} T \leq 1$.

        We know that $T$ is \ta-rigid over $A$, so \ta-rigid over $A'$ by \cref{4tauPreservedWhenPassingToQuotientAlgebras}, and $\Ext^1_{A'}(T, T) = 0$ follows.

        By \cref{4partialTiltingBound} the final thing we need to prove is that $\# T = \# A'$. Since $\# A' = \supportrank \ST = \# A - \# P$ this is seen to be equivalent to proving that $(T,P)$ is a support \ta-tilting pair. For this we consider the support $\tau^-$-tilting pair $(U, I)$ obtained from $\SF$ (using the dual of \cref{5sTauTiltFTorsBijection}). Here $U := \SI(\SF)$, and $I$ is a maximal basic injective module fulfilling $\Hom_A(\SF, I) = 0$. We now apply the bijection from \cref{5tauTauMinusBijection}.
        \[
            (U,I) \xmapsto{\quad (D(-))^\dagger \quad}(\tau^- U \oplus \nu^- I, \nu^- U_{\text{in}})
        \]
        $(\tau^- U \oplus \nu^- I, \nu^- U_{\text{in}})$ is a support \ta-tilting pair. To complete the proof we show that $T \cong \tau^- U \oplus \nu^- I$, and that $P \cong \nu^- U_\text{in}$. 

        We remind the reader of the following property of the Nakayama functor:
        \begin{equation}\label{equationNakayama}
            \Hom_A(P, T) = 0 \iff \Hom_A(T, \nu P) = 0.
        \end{equation}

        $\tau^- U \oplus \nu^- I \in \add T:$ By \cref{3extProjectivesBijection}, $\tau^- U$ is Ext-projective in $\ST$, thus a summand of $T$. Using the dual of \cref{equationNakayama}, we obtain that $\nu^- I \in {}^{\perp 0}\SF = \ST$. $\nu^- I$ is projective, so Ext-projective in $\ST$, so also a summand of $T$.

        $T \in \add (\tau^- U \oplus \nu^- I):$ Since \cref{3extProjectivesBijection} is a bijection, all non-projective summands of $T$ arise from $\tau^- U$. To check the remaining case, take an indecomposable projective $P'$ in $\ST$. Using that there are no maps from $P'$ to $\SF$ and \cref{equationNakayama}, we obtain that $\nu P' \leqslant I$. From this follows that $P' \cong \nu^- \nu P'$ is a summand of $\nu^- I$. 

        $P \in \add (\nu^- U_\text{in}):$ Using \cref{equationNakayama}, we obtain that $\Hom_A(\ST, \nu P) = 0$. Consequently, $\nu P$ is an injective module in $\SF$, thus a summand of $U_\text{in}$, and thus $P$ is a summand of $\nu^- U_\text{in}$.

        $\nu^- U_\text{in} \in \add P:$ Using the dual of \cref{equationNakayama}, we obtain that $\Hom_A(\nu^-U_\text{in}, \ST) = 0$, and thus the assertion follows.

        It follows that $\tau^- U \oplus \nu^- I \cong T$, as well as $\nu^- U_{\text{in}} \cong P$. Consequently, $(T,P)$ is a support \ta-tilting pair, and—by the above—$T$ is a tilting module over $A'$, and $\ST$ is functorially finite.
    \end{proof}
\end{theorem}

While this result is interesting in its own right, it also turns out to be a key ingredient in the construction of a Bongartz-type completion for \ta-rigid modules, and this completion is itself a key-ingredient in establishing results on \emph{mutation} of support \ta-tilting modules.

\subsection{Bongartz Completion}
\begin{proposition}\label{5bongartzHelperLemma}
    Let $\ST$ be a functorially finite torsion class, and $U$ a \ta-rigid $A$-module. Then
    \[
        U \in \add \SP(\ST) \quad \iff \quad \gen U \subseteq \ST \subseteq {}^{\perp 0}(\tau U).
    \]
    \begin{proof}
        $\implies :$ Naturally, 
        \[
            \gen U \subseteq \gen \SP(\ST) = \ST = {}^{\perp 0}\SF \subseteq {}^{\perp 0}(\tau U),
        \] where the last inclusion follows since $\tau U$ is contained in $\SF$ by \cref{3extProjectivesBijection}.

        $\impliedby :$ We know that $U \in \ST$, and that $\Hom_A(X,  \tau U) = 0$ for all $X \in \ST$. Using \cref{3AuslanderSmaloLemma}, we immediately see that $\Ext_A^1(U, \gen X) = 0$. We obtain that $U$ is Ext-projective in $\ST$, hence contained in $\add \SP(\ST)$.
    \end{proof}
\end{proposition}

\begin{lemma}
    For any \ta-rigid $A$-module $U$ we obtain a sincere functorially finite torsion class ${}^{\perp 0}(\tau U)$.
    \begin{proof}
        Since $U$ is \ta-rigid, $\tau U$ is $\tau^-$-rigid, as $\Hom_A(\tau^- \tau U, \tau U) = \Hom_A(U\dualnp, \tau U) = 0$. Dually to \cref{3tauRigidsGenerateFFTorsionClasses}, it cogenerates a functorially finite torsion-free class $\cogen \tau U$.

        Now $({}^{\perp 0}(\tau U), \cogen \tau U)$ is a torsion theory, and ${}^{\perp 0}(\tau U)$ is functorially finite by \cref{5STFfIffSFFf}. Assume that ${}^{\perp 0}(\tau U)$ is not sincere.

        In that case there exists some simple module $S$ which is not a composition factor of any module in ${}^{\perp 0}(\tau U)$. Using $I$ to refer to the indecomposable injective module whose socle equals $S$, we get that $\Hom_A({}^{\perp 0}(\tau U), I) = 0$. Consequently, $I$ is in $\cogen \tau U$. This is a contradiction since we would be able to find a (splitting) monomorphism $I \hookrightarrow (\tau U)^n$, but $\tau U$ has no injective summands.
    \end{proof}
\end{lemma}

\begin{theorem}\label{5Bongartz}
    Let $U$ be a \ta-rigid $A$-module. Then $\ST := {}^{\perp 0}(\tau U)$ is a sincere functorially finite torsion class, and $T := \SP(\ST)$ is a \ta-tilting $A$-module satisfying $U \in \add T$, and ${}^{\perp 0}(\tau T) = \gen T$.

    We call $\SP({}^{\perp 0}(\tau U))$ the Bongartz-Completion of $U$.
    \begin{proof}
        By the previous lemma $\ST$ is a sincere functorially finite torsion class. That $T$ is a \ta-tilting module follows immediately. By \cref{5bongartzHelperLemma} $U$ is in $\add T$. From this follows that 
        \[
            \gen T \subseteq {}^{\perp 0}(\tau T) \subseteq {}^{\perp 0}(\tau U) = \gen T,
        \]
        and the assertion follows.
    \end{proof}
\end{theorem}

\begin{example}
    Take some finite-dimensional $K$-algebra $A$, and a projective module $U$. Naturally, $U$ is \ta-rigid, and therefore $(\gen U, U\bong)$ is a torsion theory.

    Using Bongartz-completion, we compute $\bong(\tau U)$. Since $U$ is projective, $\tau U = 0$, and $\bong(\tau U) = A\Mod$. We obtain $\SP(\bong(\tau U)) = {}_AA$.

    At this point the reader might wonder how the process of "completing" a \ta-rigid module $U$ to a support \ta-tilting module $U \mapsto \SP(\gen U)$ interacts with the Bongartz-completion. The Bongartz-completion of $\SP(\gen U)$ and the Bongartz-completion of $U$ need not agree with each other.

    This is easy to see using the above example: The Bongartz-completion of $\SP(\gen U)$ contains $\SP(\gen U)$ as a summand, so the two cannot agree as soon as $\SP(\gen U)$ contains a non-projective summand. For example, for $A = K(1 \rightarrow 2)$ and $U = P_1 = 11$ we have $\gen U = \add\{11, 10\}$, and $\SP(\gen U) = 11 \oplus 10$.
\end{example}

The following result is an analogue of Wakamatsu's lemma. 
\begin{lemma}\label{5Wakamatsu}
    Let $T$ be \ta-rigid. Given some $A$-module $X$, write $f : T' \rightarrow X$ for an $(\add T)$-cover of $X$. 

    Then $\ker(f) \in \bong(\tau T)$.
    \begin{proof}
        Replacing $X$ with $\im f$, we can assume that $f$ is surjective. Our situation is as follows:
        \[
            \begin{tikzcd}
                0 \rar & \ker f \rar & T' \rar["f"] & X \rar & 0
            \end{tikzcd}
        \]
        We apply $\Hom_A(-, \tau T)$ to obtain 
        \[
            \begin{tikzcd}
                0 = \Hom_A(T', \tau T) \rar & \Hom_A(\ker f, \tau T) \rar & \Ext^1_A(X, \tau T) \rar["{\Ext^1(f, \tau T)}"] & \Ext^1_A(T', \tau T)
            \end{tikzcd}
        \]
        where $\Hom_A(T', \tau T) = 0$ since $T$ is \ta-rigid. Since $f$ is an $(\add T)$-cover, the induced map $(T,f) : \Hom_A(T, T') \rightarrow \Hom_A(T, X)$ is surjective. The same then applies to the induced map $(T,f) : \ul{\Hom}_A(T, T') \rightarrow \ul{\Hom}_A(T, X)$ modulo the projectives, and AR-duality gives us that $\Ext_A^1(f, \tau T)$ is injective. The result follows.
    \end{proof}
\end{lemma}

\begin{theorem}\label{5TauTiltingEquivalent}
    The following are equivalent for a \ta-rigid $A$-module $T$:
    \begin{enumerate}
        \ii $T$ is \ta-tilting.
        \ii $T$ is maximal \ta-rigid, i.e. if $T \oplus X$ is \ta-rigid for some $A$-module $X$, then $X \in \add T$.
        \ii $\bong(\tau T) = \gen T$.
        \ii If $\Hom_A(T, \tau X) = 0$, and $\Hom_A(X, \tau T) = 0$, then $X \in \add T$.
    \end{enumerate}
    \begin{proof}
        $1) \implies 2):$ If $T \oplus X$ is \ta-rigid, then $\#(T \oplus X) \leq \# A = \# T$ by \cref{4tauRigidUpperBound}, so $X \in \add T$.

        $2) \implies 3):$ Write $T'$ for the Bongartz-completion of $T$. $T$ is maximal \ta-rigid and $T \in \add T'$, so $\add T' = \add T$. Since we have that $\bong(\tau T') = \gen T'$, we are done.
        
        $3) \implies 1):$ Let $T$ be \ta-rigid with $\bong(\tau T) = \gen T$, and write $T'$ for the Bongartz-completion of $T$. We have 
        \[
            \gen T \subseteq \gen T' \subseteq \bong(\tau T') \subseteq \bong (\tau T) = \gen T,
        \]
        and hence all inclusions are equalities. Since $\gen T' = \gen T$, we obtain an $(\add T)$-cover $f : S \rightarrow T'$ such that 
        \[
            \begin{tikzcd}
                0 \rar & Y \rar & S \rar["f"] & T' \rar & 0
            \end{tikzcd}
        \] 
        is exact. By \cref{5Wakamatsu} we have that $\Hom_A(Y, \tau T) = 0$, and thus $\Hom_A(Y, \tau T') = 0$, since $\bong(\tau T) = \bong(\tau T')$. Using the Auslander-Reiten formula, we see that $\Ext_A^1(T', Y) = D\ol{\Hom}_A(Y, \tau T') = 0$, and hence the above sequence splits. In other words, $T' \in \add T$, and the result follows.

        $1) + 3) \implies 4):$ Given an $X$ such that $\Hom_A(X, \tau T) = 0$, $3)$ gives us that $X \in \gen T$. Using \cref{3AuslanderSmaloLemma} on $\Hom_A(T, \tau X) = 0$, we see that $X$ is Ext-projective in $\gen T$. By \cref{5sTauTiltFTorsBijection}, $X$ is in $\add \SP(\gen T) = \add T$. 

        $4) \implies 2):$ This is clear.
    \end{proof}
\end{theorem}

\begin{corollary}
    The following are equivalent for a \ta-rigid pair $(T, P)$ over $A$:
    \begin{enumerate}
        \ii $(T, P)$ is a support \ta-tilting pair for $A$.
        \ii If $(T \oplus X, P)$ is \ta-rigid for some $A$-module $X$, then $X \in \add T$.
        \ii $\bong(\tau T) \cap (P\bong) = \gen T$.
        \ii If $\Hom_A(T, \tau X) = 0$, and $\Hom_A(X, \tau T) = 0$, and $\Hom_A(P, X) = 0$ then $X \in \add T$.
    \end{enumerate}
    \begin{proof}
        In view of \cref{4tauPreservedWhenPassingToQuotientAlgebras}, the result follows immediately by replacing $A$ by $A/\langle e \rangle$ for an idempotent $e$ of $A$ satisfying $\add P = \add Ae$.
    \end{proof}
\end{corollary}

\section{Hasse Quiver and Mutation}

Our primary goal was to write an introductory summary to \ta-tilting theory. We focused on the backbone of the topic, and aimed to give self-contained proofs of results which are merely cited in \cite{zbMATH06293659}, as well as to sidestep the Brenner-Butler theorem. We believe that this goal has been achieved. 

In this section we will largely omit proofs. The reason for this is simple: We do not believe that we can improve on the proofs given in \cite{zbMATH06293659}. With the theory which we developed in the previous sections Section 2.3+ of \cite{zbMATH06293659} should be quite accessible to the interested reader, and not require them to look up any further sources, so we merely provide a rough overview.

In \cref{5sTauTiltFTorsBijection} we have shown that there is a bijection between support \ta-tilting modules and functorially finite torsion classes. The set of torsion classes possesses a natural structure as a partially ordered set, given by inclusion. Using this, we obtain the same poset structure on $\sttilt A$.

\begin{definition}
    We define the \emph{Hasse quiver} of $\sttilt A$ as follows:
    \begin{itemize}
        \ii The set of vertices is $\sttilt A$.
        \ii We draw an arrow from $T$ to $T'$ if $\gen T' \subsetneq \gen T$, and there exists no functorially finite torsion class $\ST$ such that $\gen T' \subsetneq \ST \subsetneq \gen T$, i.e. $\gen T'$ is maximally included in $\gen T$.
    \end{itemize}
    We have already seen Hasse quivers, mainly in \cref{3LinearA3WithRelationExample} and in \cref{3ExampleSkewedTriangle}. We may write $Q(\sttilt A)$ for the Hasse quiver of $\sttilt A$.

    Given two support \ta-tilting modules $T$, $T'$, we write $T' < T$ or $T > T'$ if $\gen T' \subsetneq \gen T$, and $T' \leq T$ if $\gen T' \subseteq \gen T$. In particular, arrows in the Hasse quiver point from larger to smaller modules.
\end{definition}

\begin{remark}\label{6RemarkTorsionLattice}
    We can apply an analogous definition and consider the Hasse quiver obtained by considering the inclusion-relationship on the set \tors $A$ of \emph{all} torsion classes. 

    It is known that $\tors A$ is a complete lattice, meaning that \emph{every} set of torsion classes has a "join", i.e. a unique least upper bound, and a "meet", i.e. unique greatest lower bound. We recommend \cite{zbMATH07431617} for a quick glimpse into the overall topic. 

    Generally speaking, the situation is not quite as simple when looking at the poset of functorially finite torsion classes, i.e. at $\sttilt A$ and its Hasse quiver $Q(\sttilt A)$. More precisely, the poset $\sttilt A$ is a complete lattice if and only if it is finite \cite{zbMATH07130859}. \Cref{6TauFiniteDefinition} is the central starting point of this topic. In the meantime, \cref{3NonFfTorsionClassesExample} provides food for thought.
\end{remark}

$Q(\sttilt A)$ possesses a bunch of desirable and interesting properties. First, we will make an intuition back from \cref{3LinearA3WithRelationExample} more precise.

\begin{definition}
    Let $A$ be a finite dimensional $K$-algebra, as usual. We call a \ta-rigid pair $(U,Q)$ over $A$ an \emph{almost complete support \ta-tilting pair} if 
    \[
        \# U + \# Q = \# A - 1.
    \]
\end{definition}

Together with \cref{6MutationTheorem3} below, which identifies the arrows of the Hasse quiver with mutations, the following result states precisely that the Hasse quiver associated with the functorially finite torsion classes of an algebra is $\#A$-regular, i.e. each vertex has exactly $\# A$ edges, and these edges are in a one-to-one correspondence with the summands of the associated support \ta-tilting pair.

\begin{theorem}
    Let $(T,P)$ be a basic support \ta-tilting pair over $A$. We decompose $T$ and $P$ into indecomposable summands: $T = \bigoplus_{k=1}^l T_k$, and $P = \bigoplus_{k=l+1}^n P_k$, where $n = \# A$.

    Then: For any indecomposable summand $T_j$ of $T$ (respectively, $P_j$ of $P$) the almost complete support \ta-tilting pair $(\bigoplus_{k \neq j} T_k,~P)$ (respectively, $(T,~\bigoplus_{k \neq j} P_k)$) can be completed to a support \ta-tilting pair different from $(T,P)$ in a unique way.
\end{theorem}

This result follows from the following theorem, which also introduces the idea of "mutation" of modules:

\begin{theorem}\label{6MutationTheorem1}
    Any basic almost complete support \ta-tilting pair $(U,Q)$ for $A$ is a direct summand of exactly two basic support \ta-tilting pairs. 

    Denoting these as $(T,P)$ and $(T', P')$, we obtain without loss of generality that 
    \[
        \gen T' = \gen U \text{ \qquad and \qquad } \gen T = \bong(\tau U) \cap Q\bong.
    \]
    In particular, we see that $T' < T$. We say that $T'$ is a \emph{left mutation} of $T$, and denote this as $T' = \mu^-_X(T)$, where $X$ is the indecomposable module such that $T = U \oplus X$.
\end{theorem}

One direction of the above statement follows easily from the previous results:

\begin{proposition}
    For every basic \ta-rigid module $T$, which is not \ta-tilting, there exist at least two basic support \ta-tilting modules which have $T$ as a direct summand.
    \begin{proof}
        By \cref{5TauTiltingEquivalent} $\gen T$ is properly contained in $\bong(\tau T)$. $\bong(\tau T)$ is a functorially finite torsion class by \cref{5Bongartz}. By \cref{5sTauTiltFTorsBijection} we obtain two different support \ta-tilting modules $\SP(\gen T)$ and $\SP(\bong(\tau T))$, which are extensions of $T$ by \cref{5bongartzHelperLemma}.
    \end{proof}
\end{proposition}

This proposition only treats pairs of the form $(T, 0)$. Given a basic \ta-rigid pair $(U, Q)$ which is not a support \ta-tilting pair, take an idempotent $e$ with $\add Q = \add Ae$. Applying the proposition over $A/\langle e \rangle$, and using \cref{5tauRigidPairsAreModules}, we obtain at least two basic support \ta-tilting pairs which have $(U,Q)$ as a direct summand.

The more difficult direction is to prove that any almost complete support \ta-tilting pair is a summand of at most two basic support \ta-tilting pairs. We omit the proof, like many others in this section, and instead merely provide a broad overview.

First, we define left and right mutation.
\begin{definition}
    Two basic support \ta-tilting pairs $(T,P)$ and $(T', P')$ for $A$ are said to be \emph{mutations} of each other if there exists a basic almost complete support \ta-tilting pair $(U,Q)$ which is a direct summand of $(T, P)$ and $(T',P')$.

    We denote this with $T' = \mu_X(T)$ where $X$ is an indecomposable $A$-module satisfying either $T = U \oplus X$, or $P = Q \oplus X$.

    Let $T = X \oplus U$, and $T'$ be support \ta-tilting $A$-modules such that $T' = \mu_X(T)$ for some indecomposable $A$-module $X$. We say that $T'$ is a

    \noindent\begin{tabular}{p{0.45\textwidth}p{0.45\textwidth}}
        left mutation of $T$ if
        \begin{enumerate}
            \ii $T > T'$
            \ii $X \notin \gen U$
            \ii $\bong(\tau U) \subseteq \bong (\tau X)$
        \end{enumerate}
        &
        right mutation of $T$ if
        \begin{enumerate}
            \ii $T < T'$
            \ii $X \in \gen U$
            \ii $\bong(\tau U) \nsubseteq \bong (\tau X)$
        \end{enumerate}
    \end{tabular}

    In each case the three listed conditions are equivalent. We know that either $T < T'$, or that $T > T'$, by \cref{6MutationTheorem1}.
\end{definition}

We note two major results:
\begin{theorem}\label{6MutationTheorem2}
    Let $U$ and $T$ be basic support \ta-tilting $A$-modules, such that $U < T$. Then:
    \begin{enumerate}
        \ii There exists a right mutation $V$ of $U$ such that $V \leq T$.
        \ii There exists a left mutation $V'$ of $T$ such that $U \leq V'$.
    \end{enumerate}
\end{theorem}

\begin{theorem}\label{6MutationTheorem3}
    Given two support \ta-tilting $A$-modules $T, T'$, the following conditions are equivalent:
    \begin{enumerate}
        \ii $T'$ is a left mutation of $T$.
        \ii $T$ is a right mutation of $T'$.
        \ii $T' < T$, and there is no support \ta-tilting module $V$, such that $T' < V < T$.
    \end{enumerate}
    In other words, the \emph{Hasse Quiver} $Q(\sttilt A)$ of the partially ordered set $\sttilt A$, and the so-called \emph{Mutation Quiver} of $A$—whose underlying vertex-set is $\sttilt A$, and which has an arrow from $T$ to $T'$ whenever $T'$ is a left mutation of $T$—coincide.
\end{theorem}

A cute little consequence of this is the following result, reminiscent of a similar result about the Auslander-Reiten quiver:
\begin{corollary}
    If $Q(\sttilt A)$ has a finite connected component $C$, then $Q(\sttilt A) = C$.
    \begin{proof}
        Fix some $U \neq A$ in $C$. Naturally, $U < A$. Therefore, we can repeatedly apply \cref{6MutationTheorem2} to obtain a chain of right mutations of support \ta-tilting modules $U = U_0 < U_1 < U_2 < \cdots \leq A$. This chain has to terminate since $C$ is finite, and therefore $A = U_i$ for some value of $i$, and $A$ is in $C$.

        Now fix some $T$ in $\sttilt A$. We repeatedly apply the same theorem to construct a chain of left mutations $A = V_0 > V_1 > V_2 > \cdots \geq T$. This chain has to terminate, so $T$ is already in $C$.
    \end{proof}
\end{corollary}

As a final remark on the theory of mutation developed in \cite{zbMATH06293659}, it is worth pointing out that $Q(\sttilt A)$ can be computed using so-called exchange sequences, which are constructed using envelopes. This turns out to be rather convenient.

\begin{theorem}
    We take a basic \ta-rigid module $U$ with $\# U = \# A - 1$, and Bongartz-complete it to a \ta-tilting module $X \oplus U$. Here $X$ is indecomposable, and $X \notin \gen U$. We define an $(\add U)$-envelope $f : X \rightarrow U'$ of $X$. Using $f$, we take the cokernel, and obtain our exchange sequence:
    \[
        \begin{tikzcd}
            X \rar["f"] & U' \rar["g"] & Y \rar & 0
        \end{tikzcd}
    \]
    Now we obtain the following:
    \begin{enumerate}
        \ii If $U$ is not sincere, then $Y = 0$, and $\mu^-_X(X \oplus U) = U$.
        \ii If $U$ is sincere, then $Y \notin \add (X \oplus U)$ is indecomposable, and $\mu^-_X(X \oplus U) = U \oplus Y$.~\footnote{That $Y$ is indecomposable is proved in \cite{zbMATH06740616}. \cite{zbMATH06293659} leaves it as a question, making for a slightly different statement.}
    \end{enumerate}
\end{theorem}

Right mutation can be computed similarly by noting that the $(-)^\dagger$-functor from \cref{5definitionDagger} acts inclusion-reversing on $Q(\sttilt A)$.

One aspect that makes the theory of mutations particularly fruitful is the sheer number of natural bijections that exist between $\sttilt A$, and related classes.

If we are given some natural bijection from $\sttilt A$ to some rather different set $G$, then we will also obtain a notion of "mutation" on $G$, allowing us to pass between objects, and find the two unique completions of "almost complete" objects. This can be very interesting, and provide deeper insights into the elements of $G$. An example of this is the so-called \emph{brick labelling} of the arrows in $Q(\sttilt A)$.

\section{Prospects}

The topic of \ta-tilting theory has garnered a substantial amount of attention in recent years. For a summary of some of these results, we recommend \cite{https://doi.org/10.48550/arxiv.2106.00426}.

In the following we will single out a few specific results which the author deemed interesting, and try to connect them to the results of the previous sections.

One particularly interesting topic are \emph{\ta-tilting finite} algebras:
\begin{definition}\label{6TauFiniteDefinition}
    We call an algebra $A$ \emph{\ta-tilting finite} if it fulfills the following equivalent conditions:
    \begin{enumerate}
        \ii There are only finitely many isomorphism classes of indecomposable \ta-rigid $A$-modules.
        \ii There are only finitely many isomorphism classes of support \ta-tilting $A$-modules, i.e. $|\sttilt A | < \infty.$
        \ii There are only finitely many isomorphism classes of basic \ta-tilting $A$-modules.
        \ii There are only finitely many torsion classes in $A\Mod$.
        \ii Every torsion class in $A\Mod$ is functorially finite.
        \ii There are only finitely many isomorphism classes of bricks in $A\Mod$.\footnote{A module $X$ is a \emph{brick} if $\End_A(X)$ is a division algebra.}
    \end{enumerate}
\end{definition}

While we called it a definition, proving that these statements are equivalent requires a significant amount of work. We will (very roughly) sketch the approach.

For a deeper understanding of this result, the reader is encouraged to look at \cite{zbMATH07130859}, as well as at \cite{zbMATH07313480}, which further expands on these results and "completes" the bijection between \emph{indecomposable} \ta-rigid modules and bricks to a bijection between $\sttilt A$ and the class of \emph{semibricks}\footnote{A \emph{semibrick} is a set $\SSS$ of isomorphism classes of bricks, such that $\Hom_A(X, Y) = 0$ for all $X \neq Y \in \SSS$.}.

\begin{proof}
    $1) \implies 2) \implies 3)$ is easy enough, and $3) \implies 1)$ is an immediate result of Bongartz-completion.

    Proving that these results are equivalent to $4)$ and $5)$ requires some heavy machinery. $4) \implies 2)$ is clear, but the rest of this result largely rests on the following theorem, which can be considered a generalisation of \cref{6MutationTheorem2} to more arbitrary torsion classes.
    \begin{theorem}[Theorem 3.1, \cite{zbMATH07130859}]
        Let $\ST$ be a functorially finite torsion class, and $\SSS$ be an arbitrary torsion class such that either $\SSS \subsetneq \ST$, or $\ST \subsetneq \SSS$. Then there exists a functorially finite torsion class $\ST'$, such that
        \begin{enumerate}
            \ii If $\SSS \subsetneq \ST$, then $\ST'$ is maximal with the property that $\SSS \subseteq \ST' \subsetneq \ST$.
            \ii If $\ST \subsetneq \SSS$, then $\ST'$ is minimal with the property that $\ST \subsetneq \ST' \subseteq \SSS$.
        \end{enumerate}
    \end{theorem}
    Now, assuming that there exists a torsion class $\SSS$ which is not functorially finite, then we are able to construct an infinite, strictly ascending chain of functorially finite torsion classes 
    \[
        0 = \ST_0 \subsetneq \ST_1 \subsetneq \ST_2 \subsetneq \cdots \subsetneq \SSS
    \]
    and this chain cannot stabilize since $\ST_i$ is functorially finite for every value of $i$, while $\SSS$ is not. This proves $4) \implies 5)$, as well as $2) \implies 5)$.

    The $5) \implies 4)$ direction follows from the observation that the intersection $\ST_\infty$ of an infinite strictly descending chain $\ST_1 \supsetneq \ST_2 \supsetneq \cdots$ of functorially finite torsion classes is a torsion class which cannot be functorially finite. Suppose otherwise. Each $\ST_i$ strictly contains $\ST_\infty$, so by the above theorem each $\ST_i$ contains a functorially finite torsion class which is minimal among those strictly containing $\ST_\infty$. By \cref{6MutationTheorem3} these minimal torsion classes correspond to the arrows of $Q(\sttilt A)$ ending in $\SP(\ST_\infty)$, so there are at most $\# A$ of them. One of them is thus contained in infinitely many $\ST_i$, hence in all of them, and hence in $\ST_\infty$, a contradiction.

    It remains to find such a chain for an algebra with infinitely many torsion classes, all of which are functorially finite by $5)$. Then $\sttilt A$ is infinite. Given any $T \in \sttilt A$, we have $T \leq A$, and repeatedly applying \cref{6MutationTheorem2} gives a chain of left mutations $A = V_0 > V_1 > V_2 > \cdots \geq T$. If this process does not terminate for some $T$, we have found an infinite strictly descending chain. Otherwise the chain ends in $T$ for every $T$, so every one of the infinitely many vertices of $Q(\sttilt A)$ can be reached by a directed path starting in $A$. Since each vertex has at most $\# A$ outgoing arrows, König's lemma yields an infinite directed path starting in $A$. Since arrows in $Q(\sttilt A)$ point from larger to smaller torsion classes, such a path gives a strictly descending chain.

    What remains is to prove the equivalence to $6)$. This follows from the following result: 
    \begin{theorem}[Theorem 4.1, \cite{zbMATH07130859}]
        There exists a bijection 
        \[
            \mbox{\rm i$\tau$-rigid}\hspace{.01in} A \rightarrow \fbrick A
        \]
        given by $X \mapsto X / \rad_{\End_A(X)}(X)$, where $\mbox{\rm i$\tau$-rigid}\hspace{.01in} A$ is the set of isomorphism classes of indecomposable \ta-rigid $A$-modules, and $\fbrick A$ is the set of isomorphism classes of bricks $S$ of $A$, such that the smallest torsion class $T(S)$ containing $S$ is functorially finite. 
    \end{theorem}
    This is enough to prove $6) \implies 1)$. To complete the equivalence it is shown that, if $S, S'$ are any bricks fulfilling $T(S) = T(S')$, then $S$ and $S'$ are isomorphic. This means in particular that, if there are only finitely many torsion classes, then there can only be finitely many bricks.
\end{proof}

This very rough outline of the results in \cite{zbMATH07130859} primarily served to illustrate three ideas which are frequently encountered in the area of \ta-tilting theory:
\begin{enumerate}
    \ii Bijections between classes of \ta-rigid modules and other classes of objects.
    \ii Mutation of support \ta-tilting modules, or generalisations of this concept.
    \ii \ta-tilting finite algebras.
\end{enumerate}

Naturally, this list is by no means exhaustive, for example, it does not even mention the important topics of silting theory or g-vectors.

We now further elaborate on \ta-tilting finite algebras, and start off with a trivial observation: Every representation-finite algebra is \ta-tilting finite. In other words, \ta-tilting finiteness can be considered a generalisation of representation-finiteness.

The immediate question is whether, perhaps, these two classes of algebras coincide. 

For many classic classes of algebras\footnote{See \emph{Chapter 10} in \cite{https://doi.org/10.48550/arxiv.2106.00426}.} the property of being representation-finite, and the property of being \ta-tilting finite coincide. Using the result that all indecomposable preprojective modules are bricks without self-extensions (see e.g. \cite[VIII. Lemma 2.7]{zbMATH02228448}), it is easy to see that a \ta-tilting finite algebra (which is not representation-finite) cannot have a preprojective component, which already rules out many common examples of algebras.

Still, it turns out that there are \ta-tilting finite algebras which are not representation-finite:

\begin{example}[\cite{https://doi.org/10.48550/arxiv.1307.6509}]\label{7WildTauFiniteAlgebrasExample}
    We present an example of an infinite class of wild \ta-tilting finite algebras.
    \[
        Q :
        \begin{tikzcd}
            &&&&& 2 \ar[dr, "\alpha"{name=a}] \\
            n \rar & n-1 \rar & n-2 \rar & \cdots \rar & 3 \ar[ur, "\beta"{name=b}] \ar[rr] && 1
            \arrow[dashed, no head, bend left, from=a, to=b, shorten >=0.8ex, shorten <=0.6ex]
        \end{tikzcd}
    \]
We define $A$ as the algebra given by the path algebra of $Q$, modulo the relation given by the composition of $\alpha$ and $\beta$. This algebra is representation-wild as long as $n \geq 9$, or in other words, about as far away from representation-finite as it could possibly be. Its global dimension equals two, but picking a sufficiently high $n$ and adding further relations makes it easily possible to create similar examples of arbitrarily high finite global dimension.

The reader is encouraged to contrast this example with \cref{3ExampleSkewedTriangle}.

We present Ringel's proof that these algebras are wild and \ta-tilting finite.
\begin{proof}
    First, let us prove that $A$ is wild. Vertex $2$ is a node in the sense of Martinez \cite[\emph{Algebras stably equivalent to l-hereditary}]{zbMATH03670383}, thus there is a natural bijection between indecomposable non-simple representations of $A$, and the indecomposable non-simple representations of the following quiver:
    \[
        Q' :
        \begin{tikzcd}
            &&&& 2' & 2'' \dar \\
            n \rar & n-1 \rar & n-2 \rar & \cdots \rar & 3 \ar[u] \ar[r] & 1
        \end{tikzcd}
    \]
    This quiver is wild for $n \geq 9$.

    Now, let us prove that $A$ is \ta-tilting finite. For this let $M$ be an indecomposable representation of $A$, such that $M_\alpha$ and $M_\beta$ are both non-zero. Write $S(2)$ for the simple module at vertex $2$. Since $M_\alpha M_\beta = 0$ and $M_\beta \neq 0$, the kernel of $M_\alpha$ is non-zero, and we obtain that $S(2)$ is a submodule of $M$. Since $M_\alpha \neq 0$, but $M_\alpha$ vanishes on the image of $M_\beta$, this image is a proper subspace of $M_2$, and we obtain that $S(2)$ is a factor module of $M$. It is now easy to construct an endomorphism $\phi$ of $M$ with image $S(2)$. It follows that $\phi^2 = 0$ as otherwise—by Fitting's lemma—$S(2)$ would be a summand of $M$, which is impossible.

    Now let $X$ be a brick over $A$. Using the above, either $X_\alpha = 0$, or $X_\beta = 0$, otherwise we can find a non-trivial nilpotent endomorphism of $X$. It follows that $X$ is either a representation of the $\mathbb{D}_n$-quiver obtained by deleting $\alpha$, or of the $\mathbb{A}_n$-quiver obtained by deleting $\beta$. Both of these quivers are well-understood and representation-finite. In other words, $A$ has only finitely many bricks, and is thus \ta-tilting finite.
\end{proof}
\end{example}

Finally, we also have the following result:

\begin{theorem}[Theorem 5.12, \cite{https://doi.org/10.48550/arxiv.1711.01785}]
    The class of \ta-tilting finite algebras is closed under taking factor algebras.
\end{theorem}

In other words, it is easy to construct rather complicated algebras—arbitrary finite global dimension, representation-wild—which are nonetheless easy to understand from a torsion-theoretic perspective.

At the time of writing, the \ta-tilting finite algebras and their classification continues to be a topic of active research.

\subsection*{Acknowledgements}
I would like to thank Yucheng Wang for carefully reading the first version of this text, and for pointing out several gaps in it. His questions led to the main corrections in this version, concerning $g$-vectors and the number of summands of tilting modules.

The remaining corrections were made with substantial help from an AI assistant (Claude, by Anthropic). It reviewed the entire text, identified several further errors and gaps—among them the false claim that $t(-)$ is always left exact, and steps missing from several proofs—and drafted the corresponding fixes. All of these changes were reviewed by the author, who takes full responsibility for the content.

\bibliography{Literatur}

\vspace{3em}

\textsc{Arne Johannsmann}

\texttt{arnejohannsmann@protonmail.ch}

\end{document}